\documentclass[12pt,a4paper,reqno,draft]{amsart}
\usepackage{amssymb}
\usepackage{latexsym}
\usepackage{exscale}

\numberwithin{equation}{section}
\newtheorem{theor}{Theorem}[section]
\newtheorem{lemma}[theor]{Lemma}
\newtheorem{defi}[theor]{Definition}
\newtheorem{corol}[theor]{Corollary}
\newtheorem{remark}[theor]{Remark}
\newtheorem{prop}[theor]{Proposition}
\newtheorem{defin}[theor]{Definition}

\newcommand{\re}{\mathbb{R}}
\newcommand{\co}{\mathbb{C}}

\newcommand{\RR}{\mathbb{R}}
\newcommand{\NN}{\mathbb{N}}
\newcommand{\HH}{\mathbb{H}}
\newcommand{\TT}{\mathbb{T}}

\newcommand{\J}{\mathcal{J}}
\newcommand{\K}{\mathcal{K}}
\newcommand{\F}{\mathcal{F}}

\newcommand{\W}{\mathcal{W}}

\newcommand{\D}{\mathcal{D}}
\newcommand{\cals}{\mathcal{S}}
\newcommand{\M}{\mathcal{M}}

\newcommand{\A}{\mathcal{A}}

\newcommand{\BMO}{{\rm BMO}}

\newcommand{\ep}{\varepsilon}
\newcommand{\normH}[1]{|\!|\!|#1|\!|\!|}
\newcommand{\bignormH}[1]{\Big|\!\Big|\!\Big|#1\Big|\!\Big|\!\Big|}

\newcommand{\bigchi}{\mathop{\mathchoice%
{\mbox{\Large$\chi$}}{\mbox{\large$\chi$}}{\mbox{\normalsize$\chi$}}%
{\mbox{\small$\chi$}}}\nolimits}

\renewcommand{\emptyset}{\mbox{\rm \O}}
\renewcommand{\Re}{{\rm Re}\,}

\def\div{\mathop{\rm div}}
\def\Int{\mathop{\rm Int}}
\def\supp{\mathop{\rm supp}}

\newcommand{\expt}[1]{e^{\textstyle #1}}

\newcommand{\aver}[1]{-\hskip-0.46cm\int_{#1}}

\newcommand{\off}[2]{\mathcal{O}\big(L^{#1}-L^{#2}\big)}
\newcommand{\offw}[2]{\mathcal{O}\big(L^{#1}(w)-L^{#2}(w)\big)}
\newcommand{\fullx}[2]{\mathcal{F}\big(L^{#1}-L^{#2}\big)}

\newcommand{\dec}[1]{\Upsilon\!\left(#1\right)}

\begin{document}
\allowdisplaybreaks

\title[Weighted norm inequalities and elliptic
operators]{Weighted norm inequalities,  off-diagonal estimates and
elliptic operators
\\[.2cm]
{\footnotesize Part III:  Harmonic analysis of  elliptic operators}}

\author{Pascal Auscher}

\address{Pascal Auscher
\\
Universit\'e de Paris-Sud et CNRS UMR 8628
\\
91405 Orsay Cedex, France} \email{pascal.auscher@math.p-sud.fr}

\author{Jos\'e Mar{\'\i}a Martell}

\address{Jos\'e Mar{\'\i}a Martell
\\
Instituto de Matem\'aticas y F{\'\i}sica Fundamental
\\
Consejo Superior de Investigaciones Cient{\'\i}ficas
\\
C/ Serrano 123
\\
28006 Madrid, Spain}

\address{\null\vskip-.7cm and\vskip-.7cm\null}

\address{Departamento de Matem\'aticas \\ Universidad Aut\'onoma de Madrid \\
28049 Madrid, Spain } \email{chema.martell@uam.es}

\thanks{This work was partially supported by the European Union
(IHP Network ``Harmonic Analysis and Related Problems'' 2002-2006,
Contract HPRN-CT-2001-00273-HARP). The second author was also
supported by MEC ``Programa Ram\'on y Cajal, 2005'' and by MEC Grant
MTM2004-00678.}

\date{\today}
\subjclass[2000]{42B20, 42B25, 47A60}

\keywords{Muckenhoupt weights, elliptic operators in divergence form,
singular non-integral operators, holomorphic functional calculi,
square functions, square roots of elliptic operators, Riesz
transforms, maximal regularity, commutators with bounded mean
oscillation functions.}

\begin{abstract}
This is the third part of a series of four  articles on weighted norm
inequalities, off-diagonal estimates and elliptic operators. For $L$
in some class of elliptic operators, we study
 weighted norm $L^p$ inequalities  for  singular
``non-integral'' operators arising from $L$ ; those are the operators
$\varphi(L)$ for bounded holomorphic functions $\varphi$, the Riesz
transforms $\nabla L^{-1/2}$ (or $(-\Delta)^{1/2}L^{-1/2}$) and its
inverse $L^{1/2}(-\Delta)^{-1/2}$, some  quadratic functionals
$g_{L}$ and $G_{L}$ of  Littlewood-Paley-Stein type and also some
vector-valued inequalities such as the ones involved for maximal
$L^p$-regularity. For each, we obtain sharp or nearly sharp ranges of
$p$ using the general theory  for boundedness of  Part I and the
off-diagonal estimates of Part II.   We also obtain commutator
results with BMO functions.
\end{abstract}

\maketitle

\null\vskip-1.7cm\null

\begin{quote}
{\footnotesize\tableofcontents}
\end{quote}

\section{Introduction}

In this part, we consider divergence form uniformly elliptic complex operators $ L
=
-\div(A\,\nabla )
$ in $\RR^n$ and  we are interested in weighted  $L^p$ estimates
for:
\begin{list}{$(\theenumi)$}{\usecounter{enumi}\leftmargin=.8cm
\labelwidth=0.7cm\itemsep=0.2cm\topsep=.2cm
\renewcommand{\theenumi}{\alph{enumi}}}

\item $\varphi(L)$ with  bounded holomorphic functions $\varphi$ on
sectors (Section \ref{sec:fc}).

\item The  square root $L^{1/2}$ compared to the ones for $\nabla$
and, in particular, the Riesz transforms $\nabla L^{-1/2}$ (Sections
\ref{sec:RT}, \ref{sec:RRT}).

\item Typical square functions ``\`a la'' Littlewood-Paley-Stein: one, $g_{L}$, using
only  functions of $L$,  and  the other, $G_{L}$, combining functions
of $L$ and the gradient operator (Section \ref{sec:sf}).

\item Vector-valued inequalities for the operators above and  the so-called $R$-bounded\-ness of the analytic semigroup $\{e^{-z\, L} \}$ which is linked to maximal regularity (Section \ref{sec:vv}).
\end{list}

Let us stress that those operators  may not be representable with
``usable'' kernels: they are ``non-integral''.  But they  still are
singular in the sense that they are of order 0.  Hence, usual methods
for singular integrals have to be strengthened.  The unweighted $L^p$
estimates  are described in \cite{Aus} for the operators in
$(a)-(c)$, with emphasis on the sharpness of the ranges of $p$. The
instrumental tools are two criteria for $L^p$ boundedness, valid in
spaces of homogeneous type: one was a sharper and simpler version of
a theorem by Blunck and Kunstmann \cite{BK1} in the spirit of
H\"ormander's criterion via the Calder\'on-Zygmund decomposition, and
the other one a criterion of the first author, Coulhon, Duong and
Hofmann \cite{ACDH} in the spirit of Fefferman and Stein's sharp
maximal function via a good-$\lambda$ inequality. The main interest
of those results were that they yield $L^p$ boundedness of
(sub)linear operators on   spaces of homogeneous type  for $p$ in an
arbitrary interval.  Such theorems  are  extended in Part I of our
series \cite{AM1} to obtain weighted $L^p$ bounds for the operator
itself, its commutators with a BMO function and also vector-valued
expressions.

In Part II \cite{AM2}, we studied one-parameter families of operators
satisfying   local $L^p-L^q$ estimates called off-diagonal estimates
on balls (the setting is that of a space of homogeneous type). Among
other things,  such estimates  imply  uniform $L^p$-boundedness and
are stable under composition.  In case of one-parameter semigroups,
we showed that as soon as there exists \emph{one} pair $( p,q)$ of
indices with $p<q$ for which these local $L^p-L^q$ estimates hold,
then they hold for \emph{all} pairs of indices taken in the interior
of the range of $L^p$ boundedness. This fact is of utmost importance
for applications as we often need to play with exponents. We showed
that such estimates pass from  the unweighted case to the weighted
case. Eventually, we made a thorough study of weighted off-diagonal
estimates on balls for  the semigroup arising from the operator $ L$
above.

\

Our strategy here has the same   two steps  in each of the four
situations described above. The first step consists in obtaining a
first range of exponents $p$ (depending on the weight) by applying
the abstract machinery from Part I. This range turns out to be the
best possible for both classes of operators and weights.

However, given one operator and one weight, the range of $p$ obtained
above may not be sharp, and this leads us to the second step. The
sharp range is in fact related to the one for weighted off-diagonal
estimates established in Part II. At this point, we use  the main
results of Part I in the Euclidean space but now equipped with the
doubling measure $w(x)\, dx$.

We wish to point out that some of our results can be obtained by
different methods  (essentially from geometric theory of Banach
spaces) once the bounded holomorphic functional calculus is
established in $(a)$. We give the references in the text.

We wish to say that our proofs are technically  simpler than the ones
in \cite{Aus} even for the unweighted case, because the notion of
off-diagonal estimates used here is more appropriate.

Finally, thanks to the general results in Part I,  the same
technology allows us to prove in passing  weighted $L^p$ estimates
for commutators of the operators in $(a)-(c)$ with BMO functions in
the same ranges of exponents (see Section \ref{sec:commutators}).

\section{General criteria for boundedness and the set $\W_w(p_0,q_0)$}\label{sec:general}

The underlying space is the Euclidean setting  $\RR^n$ equipped with
Lebesgue measure or a doubling measure obtained from an $A_{\infty}$
weight. We state   two results used in this work, referring  to \cite{AM1} for  statements  in stronger form and for references to earlier works.

Given a ball $B$, we write
$$
\aver{B} h\,dx
=
\frac1{|B|}\,\int_B h(x)\,dx.
$$
Let  us introduce some classical classes of weights. Let $w$ be a
weight (that is a non negative locally integrable function) on
$\RR^n$. We say that $w\in A_p$, $1<p<\infty$, if there exists a
constant $C$ such that for every ball $B\subset\re^n$,
$$
\Big(\aver{B} w\,dx\Big)\,
\Big(\aver{B} w^{1-p'}\,dx\Big)^{p-1}\le C.
$$
For $p=1$, we say that $w\in A_1$ if there is a constant $C$ such
that for every ball $B\subset \re^n$,
$$
\aver{B} w\,dx
\le
C\, w(y),
\qquad \mbox{for a.e. }y\in B.
$$
The reverse H\"older classes are defined in the following way: $w\in
RH_{q}$, $1< q<\infty$, if there is a constant $C$ such that  for any
ball $B$,
$$
\Big(\aver{B} w^q\,dx\Big)^{\frac1q}
\le C\, \aver{B} w\,dx.
$$
The endpoint $q=\infty$ is given by the condition $w\in RH_{\infty}$
whenever there is a constant $C$ such that for any ball $B$,
$$
w(y)\le C\, \aver{B} w\,dx,
\qquad \mbox{for a.e. }y\in B.
$$

The following facts are  well-known (see for instance \cite{GR, Gra}).
\begin{prop}\label{prop:weights}\
\begin{enumerate}
\renewcommand{\theenumi}{\roman{enumi}}
\renewcommand{\labelenumi}{$(\theenumi)$}
\addtolength{\itemsep}{0.2cm}

\item $A_1\subset A_p\subset A_q$ for $1\le p\le q<\infty$.

\item $RH_{\infty}\subset RH_q\subset RH_p$ for $1<p\le q\le \infty$.

\item If $w\in A_p$, $1<p<\infty$, then there exists $1<q<p$ such
that $w\in A_q$.

\item If $w\in RH_q$, $1<q<\infty$, then there exists $q<p<\infty$ such
that $w\in RH_p$.

\item $\displaystyle A_\infty=\bigcup_{1\le p<\infty} A_p=\bigcup_{1<q\le
\infty} RH_q $

\item If $1<p<\infty$, $w\in A_p$ if and only if $w^{1-p'}\in
A_{p'}$.

\item If $w\in A_{\infty}$, then the measure $dw=w\, dx$ is a Borel doubling measure.

\end{enumerate}
\end{prop}

If the Lebesgue measure is replaced by a Borel doubling measure
$\mu$, then all the above properties remain valid with the notation
change
\cite{ST}.

Given $1\le p_0<q_0\le \infty$ and $w\in A_\infty$ (with respect to a Borel doubling measure $\mu$) we define the set
$$
\W_w(p_0,q_0)
=
\big\{
p: p_0<p<q_0, w\in A_{\frac{p}{p_0}}\cap
RH_{\left(\frac{q_0}{p}\right)'}
\big\}.
$$
If $w=1$, then $\W_{1}(p_0,q_0)= (p_0,q_0)$. As it is shown in \cite{AM1}, if not empty, we have
$$\W_w(p_0,q_0)=\Big(p_0\,r_w,
\frac{q_0}{(s_w)'}\Big)
$$
where
$$
r_w=\inf\{r\ge 1\, : \, w\in A_r\},
\qquad\qquad
s_w=\sup\{s>1\, : \, w\in RH_s\}.
$$

We    use the following notation: if $B$ is a ball
with radius $r(B)$ and $\lambda>0$, $\lambda\, B$ denotes the
concentric ball with radius $r(\lambda \, B)= \lambda\, r(B)$,
$C_{j}(B)=2^{j+1}\, B\setminus 2^j\, B$ when $j\ge 2$, $C_{1}(B)=4B$,
and
$$
\aver{C_j(B)} h\,d\mu
=
\frac1{\mu(2^{j+1}B)}\,\int_{C_{j}(B)} h\,d\mu.
$$

\begin{theor} \label{theor:main-w} Let $\mu$ be a doubling Borel measure on $\re^n$ and $1\le p_0<q_0\le \infty$. Let $T$ be a sublinear
operator acting on $L^{p_0}(\mu)$, $\{\A_r\}_{r>0}$ a
family of operators acting from a subspace $\D$ of $L^{p_0}(\mu)$ into $L^{p_{0}}(\mu)$ and   $S$  an operator from $\D$ into the space of measurable functions on $\re^n$.  Assume that
\begin{equation}\label{T:I-A}
\Big(\aver{B}
|T(I-\A_{r(B)})f|^{p_0}\,d\mu\Big)^{\frac1{p_0}}
\le
\sum_{j\ge 1} g(j)\,
\Big( \aver{2^{j+1}\,B}
|S f|^{p_0}\,d\mu\Big)^{\frac1{p_0}},
\end{equation}
and
\begin{equation}\label{T:A}
\Big(\aver{B} |T\A_{r(B)}f|^{q_0}\,d\mu\Big)^{\frac1{q_0}}
\le
\sum_{j\ge 1} g(j)\,
\Big( \aver{2^{j+1}\,B}
|T f|^{p_0}\,d\mu\Big)^{\frac1{p_0}},
\end{equation}
for all $f\in \D$, all ball $B$ where $r(B)$ denotes its radius for
some $g(j)$ with $\sum g(j)<\infty$ \textup{(}with usual changes if
$q_{0}=\infty$\textup{)}. Let $p\in \W_w(p_0,q_0)$, that is,
$p_0<p<q_0$ and $w\in A_{\frac{p}{p_0}}\cap
RH_{\left(\frac{q_0}{p}\right)'}$. There is a  constant $C$ such that
for all  $f\in  \D$
\begin{equation}\label{est:main-th} \|T f\|_{L^p(w)}
\le
C\, \|S f\|_{L^p(w)}.
\end{equation}
\end{theor}

An operator acting from $A$ to $B$ is just a map from $A$ to $B$.
Sublinearity means $\vert T(f+g) \vert \le \vert Tf\vert + \vert
Tg\vert$ and $\vert T(\lambda f)\vert = \vert \lambda\vert \vert
T(f)\vert$ for all $f,g$ and $\lambda\in \re$ or $\co$ (although the
second property is not needed in this section). Next, $L^p(w)$ is the
space of complex valued functions in  $ L^p(dw)$  with  $dw=w\,d\mu$.
However, all this extends to   functions  valued in a Banach space.

\begin{remark}\rm
In the applications below, we have, either $S f=f$ with $f\in
L^\infty_{c}$ the space of compactly supported bounded functions on
$\RR^n$, or $S f=\nabla f$ with $f\in \cals$ the Schwartz class on
$\RR^n$ (see Section \ref{sec:RRT}).
\end{remark}

Let us recall that  the doubling order $D$ of a doubling measure $\mu$ is the smallest number $\kappa\ge 0$ such that there exists  $C\ge 0$ for which
$\mu(\lambda\,B)\le C_\mu\,\lambda^\kappa\,\mu(B)$ for every ball $B$ and
for any $\lambda>1$.

The other criterion we are going to use is the following.

\begin{theor}\label{theor:SHT:small}
Let $\mu$ be a doubling Borel measure on $\re^n$,  $D$ its doubling order
and $1\le p_0< q_0\le \infty$. Suppose that $T$ is a sublinear
operator bounded on $L^{q_0}(\mu)$ and that $\{\A_r\}_{r>0}$ is
family of linear operators acting from $L^\infty_{c}$ into $L^{q_0}(\mu)$.  Assume that for $j\ge 2$,
\begin{equation}\label{SHT:small:T:I-A}
\Big( \aver{C_j(B)} |T(I-\A_{r(B)})f|^{p_0}\,d\mu\Big)^{\frac1{p_0}}
\le
g(j)\,\Big(\aver{B} |f|^{p_0}\,d\mu\Big)^{\frac1{p_0}}
\end{equation}
and for $j\ge 1$,
\begin{equation}\label{SHT:small:A}
\Big(
\aver{C_j(B)} |\A_{r(B)}f|^{q_0}\,d\mu\Big)^{\frac1{q_0}}
\le
g(j)\,\Big(\aver{B} |f|^{p_0}\,d\mu\Big)^{\frac1{p_0}}
\end{equation}
for all ball $B$ with $r(B)$ its radius and for all $f\in L^\infty_{c}$ supported in
$B$. If $\sum_j g(j)\,2^{D\,j}<\infty$ then $T$ is of weak type
$(p_0,p_0)$ and hence $T$ is of strong type $(p,p)$  for all
$p_0<p<q_0$. More precisely,  there exists a constant $C$ such that for all $f \in L^\infty_{c}$
$$
\|Tf\|_{L^p(\mu)} \le C\, \| f\|_{L^p(\mu)}.
$$
\end{theor}

Again, the statement has a vector-valued extension  for linear operators  acting on and into  $L^p$ functions valued in a Banach space.

\begin{remark}\rm  Notice the symmetry between   \eqref{T:I-A} and \eqref{SHT:small:T:I-A}.
\end{remark}

\section{Off-diagonal estimates}\label{sec:off}

We first introduce the class of elliptic operators considered in this work.
Let $A$ be an $n\times n$ matrix of complex and
$L^\infty$-valued coefficients defined on $\re^n$. We assume that
this matrix satisfies the following ellipticity (or \lq\lq
accretivity\rq\rq) condition: there exist
$0<\lambda\le\Lambda<\infty$ such that
$$
\lambda\,|\xi|^2
\le
\Re A(x)\,\xi\cdot\bar{\xi}
\quad\qquad\mbox{and}\qquad\quad
|A(x)\,\xi\cdot \bar{\zeta}|
\le
\Lambda\,|\xi|\,|\zeta|,
$$
for all $\xi,\zeta\in\co^n$ and almost every $x\in \re^n$. We have used the notation
$\xi\cdot\zeta=\xi_1\,\zeta_1+\cdots+\xi_n\,\zeta_n$ and therefore
$\xi\cdot\bar{\zeta}$ is the usual inner product in $\co^n$. Note
that then
$A(x)\,\xi\cdot\bar{\zeta}=\sum_{j,k}a_{j,k}(x)\,\xi_k\,\bar{\zeta_j}$.
Associated with this matrix we define the second order divergence
form operator
$$
L f
=
-\div(A\,\nabla f),
$$
which is understood in the standard weak sense as a maximal-accretive operator on $L^2(\RR^n,dx)$ with domain $\D(L)$ by means of a
sesquilinear form.

The operator $-L$ generates a $C^0$-semigroup
$\{e^{-t\,L}\}_{t>0}$ of contractions on $L^2(\RR^n,dx)$.
Define  $\vartheta\in[0,\pi/2)$ by,
$$
\vartheta = \sup\{ \big|\arg \langle  Lf,f\rangle\big|\, : \, f\in\mathcal{D}(L)\}.
$$
Then the semigroup   has an analytic extension to a complex semigroup
$\{e^{-z\,L}\}_{z\in\Sigma_{\pi/2- \vartheta}}$ of contractions on
$L^2(\RR^n,dx)$. Here we have written for $0<\theta<\pi$,
$$
\Sigma_{\theta}
=
\{z\in\co^*:|\arg z|<\theta\}.
$$

Let $w\in A_{\infty}$. Here and thereafter, we write $L^p(w)$ for
$L^p(\RR^n, wdx)$ and if $w=1$, we drop $w$ in the notation. We
define  $\widetilde \J_{w}(L)$ and $\widetilde \K_{w}(L)$ as the
(possibly empty) intervals of those exponents $p\in [1,\infty]$ such
that $\{e^{-t\,L}\}_{t>0}$ is a bounded set in  $\mathcal{L}(L^p(w))$
and $\{\sqrt t\, \nabla e^{-t\,L}\}_{t>0}$ is a bounded set in
$\mathcal{L}(L^p(w))$ respectively (where $\mathcal{L}(X)$ is the
space of linear continuous maps on a Banach space $X$).

We  extract from \cite{Aus, AM2}  some definitions and results (sometimes in weaker form) on  unweighted  and weighted off-diagonal estimates. See there for details and more precise statements.  Set $d(E,F)=\inf \{|x-y|\, : \, x\in E, y \in F\}$ where $E, F$ are subsets of $\RR^n$.

\begin{defin}\label{def:full}
Let $1\le p\le q \le \infty$. We say that a family $\{T_t\}_{t>0}$ of
sublinear operators satisfies  $L^p-L^q$  {full off-diagonal
estimates}, in short $T_{t}\in \fullx{p}{q}$,   if for some $c>0$,
for all closed sets $E$ and $F$, all $f$ and all $t>0$ we
have\footnote[2]{Here and thereafter, for two positive quantities $A,
B$, by $A\lesssim B$ we mean that there exists a constant $C> 0$
(independent of the various parameters) such that $A\le C B$.}
\begin{equation}\label{eq:offLpLq}
\Big(\int_{F}|T_t (\bigchi_{E}\,  f)|^q\, dx\Big)^{\frac 1q} \lesssim  t^{- \frac 1 2 (\frac n p - \frac n q )}
\expt{-\frac{c\, d^2(E,F)}{t}} \Big( \int_{E}|f|^p\, dx\Big)^{\frac 1 p}.
\end{equation}
\end{defin}

We set  $\dec{s}=\max\{s,s^{-1}\}$ for $s>0$.  Given a ball $B$, recall that $C_j(B)=2^{j+1}\,B\setminus
2^j\,B$ for $j\ge 2$   and if $w\in A_{\infty}$ we use the notation
$$
\aver{B} h\,dw
=
\frac1{w(B)}\,\int_B h\,dw,
\aver{C_j(B)} h\,dw
=
\frac1{w(2^{j+1}B)}\,\int_{C_{j}(B)} h\,dw.
$$

\begin{defi}\label{defi:off-d:weights}
Given $1\le p\le q\le \infty$ and any weight $w\in A_{\infty}$, we say that a
family of sublinear operators $\{T_t\}_{t>0}$ satisfies
$L^{p}(w)-L^{q}(w)$ off-diagonal estimates on balls, in short $T_t \in\offw{p}{q}$, if there
exist $\theta_1, \theta_2>0$ and $c>0$ such that for every $t>0$ and
for any ball $B$ with radius $r$ and all $f$,
\begin{equation}\label{w:off:B-B}
\Big(\aver{B} |T_t( \bigchi_B\, f) |^{q}\,dw\Big)^{\frac 1q}
\lesssim
\dec{\frac{r}{\sqrt{t}}}^{\theta_2} \,\Big(\aver{B}
|f|^{p}\,dw\Big)^{\frac 1p};
\end{equation}
and, for all $j\ge 2$,
\begin{equation}\label{w:off:C-B}
\Big(\aver{B}|T_t( \bigchi_{C_j(B) }\, f) |^{q}\,dw\Big)^{\frac 1q}
\lesssim
2^{j\,\theta_1}\, \dec{\frac{2^j\,r}{\sqrt{t}}}^{\theta_2}\,
\expt{-\frac{c\,4^{j}\,r^2}{t}} \,
\Big(\aver{C_j(B)}|f|^{p}\,dw\Big)^{\frac 1p}
\end{equation}
and
\begin{equation}\label{w:off:B-C}
\Big(\aver{C_j(B)}|T_t( \bigchi_B\, f) |^{q}\,dw\Big)^{\frac 1q}
\lesssim
2^{j\,\theta_1}\, \dec{\frac{2^j\,r}{\sqrt{t}}}^{\theta_2}\,
\expt{-\frac{c\,4^{j}\,r^2}{t}}
\,\Big(\aver{B}|f|^{p}\,dw\Big)^{\frac 1p}.
\end{equation}
\end{defi}
Let us make some relevant comments for this work (see \cite{AM2} for
further details).

\begin{list}{$\bullet$}{\leftmargin=.6cm
\labelwidth=0.7cm\itemsep=0.3cm\topsep=.3cm}

\item In the Gaussian factors the value of $c$ is irrelevant as long as it remains non negative. We will
freely use the same letter from line to line even if its value changes.

\item These definitions extend to complex families $\{T_z\}_{z\in \Sigma_{\theta}}$
with $t$ replaced by $|z|$ in the estimates.

\item In both definitions,
$T_{t}$ may only be defined on a dense subspace $\D$ of $L^p$ or
$L^p(w)$ ($1\le p<\infty$) that is stable by truncation by indicator
functions of measurable sets (for example, $L^p\cap L^2$, $L^p(w)
\cap L^2$ or $L^\infty_{c}$).

\item If $q=\infty$,
one should adapt the definitions in the usual straightforward way.

\item $L^1(w)-L^\infty(w)$  off-diagonal estimates on balls
are equivalent to pointwise  Gaussian upper bounds  for the
kernels of $T_{t}$.

\item Both notions are stable by composition:
$T_t\in\offw{q}{r}$ and $S_t\in\offw{p}{q}$ then $T_t\circ S_t\in
\offw{p}{r}$ when $1\le p
\le q \le r\le\infty$ and similarly for full off-diagonal estimates.

\item When $w=1$, $L^p-L^q$ off-diagonal estimates on balls  are equivalent to $L^p-L^q$ full off-diagonal estimates.

\item Notice the symmetry between \eqref{w:off:C-B} and
\eqref{w:off:B-C}.

\end{list}

If $I$ is a subinterval of $[1,\infty]$, $\Int I$ denotes the interior in $\RR$ of
$I\cap \RR$.

\begin{prop} \label{prop:sgfull}
Fix $m\in \NN$ and $0<\mu <\pi/2-\vartheta$.
\begin{list}{$(\theenumi)$}{\usecounter{enumi}\leftmargin=.8cm
\labelwidth=0.7cm\itemsep=0.3cm\topsep=.3cm
\renewcommand{\theenumi}{\alph{enumi}}}

\item   There exists a non empty maximal  interval in $[1,\infty]$,  denoted by $\J(L)$,
such that if $p, q \in
\J(L)$ with $p\le q$, then
$\{(zL)^me^{-z\,L}\}_{z\in \Sigma_{\mu}}$ satisfies $L^p-L^q$
full off-diagonal estimates  and is a  bounded set in  $\mathcal{L}(L^p)$.
Furthermore, $\J(L)\subset \widetilde \J(L)$ and $\Int\J(L)=\Int\widetilde\J(L)$.

\item There exists a non empty maximal interval of $[1,\infty]$, denoted by $\K(L)$,  such that
if $p, q \in
\K(L)$ with $p\le q$, then
$\{ \sqrt z \, \nabla (zL)^me^{-z\,L}\}_{z\in \Sigma_{\mu}}$ satisfies $L^p-L^q$
full off-diagonal estimates and  is a  bounded set in  $\mathcal{L}(L^p)$.
  Furthermore, $\K(L) \subset \widetilde \K(L)$ and $\Int\K(L)=\Int\widetilde\K(L)$.

\item $\K(L) \subset \J(L)$ and, for $p<2$, we have $p\in \K(L)$ if and only if $p\in \J(L)$.

\item Denote by $p_{-}(L), p_{+}(L)$ the lower and upper bounds of \,$\J(L)$ \textup{(}hence,
of \/ $\Int\widetilde\J(L)$ also\textup{)} and by $q_{-}(L),
q_{+}(L)$ those of \,$\K(L)$ \textup{(}hence, of
\,$\Int\widetilde\K(L)$ also\textup{)}. We have $p_{-}(L)=q_{-}(L)$
and $(q_{-}(L))^*\le p_{+}(L)$.

\item If $n=1$,  $\J(L) =\K(L)=[1,\infty]$.

\item If $n=2$, $\J(L)=[1,\infty]$ and $\K(L)\supset [1,q_{+}(L))$ with $q_{+}(L)>2$.

\item If $n\ge 3$,  $p_{-}(L)< \frac{2n}{n+2}$,  $p_{+}(L)> \frac{2n}{n-2}$ and $q_{+}(L)>2$.
\end{list}
\end{prop}

We have set $q^*= \frac{q\, n }{n - q}$, the Sobolev exponent of $q$ when $q<n$ and $q^*=\infty$ otherwise.

\begin{prop} \label{prop:sg-w:extension}
Fix $m\in \NN$ and $0<\mu <\pi/2-\vartheta$.  Let $w \in A_{\infty}$.
\begin{list}{$(\theenumi)$}{\usecounter{enumi}\leftmargin=.8cm
\labelwidth=0.7cm\itemsep=0.3cm\topsep=.3cm
\renewcommand{\theenumi}{\alph{enumi}}}
\item  Assume $\W_{w}\big(p_{-}(L), p_{+}(L)\big)\ne \emptyset$. There is a maximal interval of $[1,\infty]$,  denoted by $\J_w(L)$,   containing $\W_{w}\big(p_{-}(L), p_{+}(L)\big)$, such that  if $p, q \in
\J_{w}(L)$ with $p\le q$, then
$\{(zL)^me^{-z\,L}\}_{z\in \Sigma_{\mu}}$ satisfies $L^p(w)-L^q(w)$
off-diagonal estimates on balls and is a  bounded set in  $\mathcal {L} (L^p(w))$.
  Furthermore, $\J_{w}(L)\subset \widetilde\J_{w}(L)$ and $\Int\J_{w}(L)=\Int\widetilde\J_{w}(L)$.

\item Assume $\W_{w}\big(q_{-}(L), q_{+}(L)\big)\ne \emptyset$. There exists a maximal interval of $[1,\infty]$, denoted by $\K_w(L)$, containing $\W_{w}\big(q_{-}(L), q_{+}(L)\big)$ such that if $p, q \in
\K_{w}(L)$ with $p\le q$, then
$\{ \sqrt z \, \nabla (zL)^me^{-z\,L}\}_{z\in \Sigma_{\mu}}$ satisfies $L^p(w)-L^q(w)$
off-diagonal estimates on balls and is a bounded set in $\mathcal{L}(L^p(w))$.  Furthermore, $\K_{w}(L)\subset \widetilde \K_{w}(L)$ and $\Int\K_{w}(L)=\Int\widetilde\K_{w}(L)$.

\item  Let $n\ge 2$. Assume $ \W_{w}\big(q_{-}(L), q_{+}(L)\big)\ne \emptyset$. Then $\K_{w}(L)\subset \J_{w}(L)$. Moreover,   $\inf \J_{w}(L) = \inf \K_{w}(L)$ and
$(\sup \K_{w}(L))^*_w \le  \sup \J_{w}(L) $.
\item If $n=1$,  the intervals $\J_{w}(L)$ and $\K_{w}(L)$ are the same and contain $(r_{w},\infty]$ if $w\notin A_{1}$ and are equal to  $[1,\infty]$ if $w\in A_{1}$.
\end{list}
\end{prop}

We have set $q^*_w= \frac{q\, n \, r_{w}}{n\, r_{w} - q}$ when $q<n\, r_{w}$ and $q^*_{w}=\infty$ otherwise. Recall that $r_{w}=\inf\{r\ge 1\, : \, w\in A_{r}\}$ and
also that $s_{w}=\sup\{s>1\, : \, w\in RH_{s}\}$.

Note that  $\W_{w}\big(p_{-}(L), p_{+}(L)\big)\ne \emptyset$ means
$\frac{p_{+}(L)}{p_{-}(L)}  >     r_{w}(s_{w})'$.  This is a
compatibility condition between $L$ and $w$. Similarly,
$\W_{w}\big(q_{-}(L), q_{+}(L)\big)\ne \emptyset $ means
$\frac{q_{+}(L)}{q_{-}(L)}  >  r_{w}(s_{w})'$, which is a more
restrictive condition on $w$ since $q_{-}(L)=p_{-}(L)$ and
$q_{+}(L)\le p_{+} (L)$.

In the case of real operators, $\J(L)=[1,\infty]$ in all dimensions
because the kernel $e^{-t\, L}$ satisfies a pointwise Gaussian upper
bound.  Hence $\W_{w}\big(p_{-}(L), p_{+}(L)\big)= (r_{w},\infty)$.
If $w\in A_{1}$, then one has that $\J_{w}(L)=[1,\infty]$. If
$w\notin A_{1}$, since the kernel is also positive and satisfies a
similar pointwise lower bound, one has $\J_{w}(L) \subset
(r_{w},\infty] $.  Hence, $\Int \J_{w}(L)= \W_{w}\big(p_{-}(L),
p_{+}(L)\big)$.

The situation may change for complex operators. But we lack of examples to say whether or not
  $\J_{w}(L)$ and $ \W_{w}\big(p_{-}(L), p_{+}(L)\big)$  have  different endpoints.

\begin{remark}\label{remark:inf-gen}
\rm
Note that by density of $L^\infty_{c}$ in the spaces $L^p(w)$ for
$1\le p <\infty$, the various extensions of $e^{-z\, L}$ and $\nabla
e^{-z\, L}$ are all consistent. We keep the above notation to denote
any such extension.  Also, we showed in \cite{AM2} that as long as
$p\in \J_{w}(L)$ with $p\ne\infty$, $\{e^{-t\, L}\}_{t>0}$ is strongly
continuous on $L^p(w)$, hence it has an
infinitesimal generator in $L^p(w)$, which is of type $\vartheta$.
\end{remark}

\textit{From now on, $L$ denotes an operator as defined in this
section with the four numbers $p_{-}(L)=q_{-}(L)$ and $p_{+}(L)$,
$q_{+}(L)$.  We  often drop $L$ in the notation: $p_{-}=p_{-}(L)$, $
\dots$. For a given weight $w\in A_{\infty}$, we set
$\W_w\big(p_-,p_+\big)=(\widetilde{p}_-,\widetilde{p}_+)$
\textup{(}when it is not empty\textup{)}  and $\Int
\J_w(L)=(\widehat{p}_-,\widehat{p}_+)$. We have $\widehat{p}_-\le
\widetilde{p}_-<\widetilde{p}_+\le \widehat{p}_+$. Similarly,  we set
$\W_w\big(q_-,q_+\big)=(\widetilde{q}_-,\widetilde{q}_+)$
\textup{(}when it is not empty\textup{)} and
$\Int\K_w(L)=(\widehat{q}_-,\widehat{q}_+)$. We have
$\widehat{q}_-\le \widetilde{q}_-<\widetilde{q}_+\le \widehat{q}_+$.
}

\section{Functional Calculi}\label{sec:fc}

Let  $\mu\in(\vartheta,\pi)$  (do not confuse with the measure $\mu$
used  in Section \ref{sec:general}) and
 $\varphi$ be a holomorphic function in $\Sigma_{\mu}$  with the following decay
\begin{equation}\label{decay:varphi}
|\varphi(z)|
\le
c\,|z|^s\,(1+|z|)^{-2\,s},
\qquad
z\in\Sigma_\mu,
\end{equation}
for some $c$, $s>0$. Assume that $\vartheta<\theta<\nu<\mu<\pi/2$. Then we have
\begin{equation}\label{phi-L}
\varphi(L)
=
\int_{\Gamma_+} e^{-z\,L}\,\eta_+(z)\,dz +\int_{\Gamma_-}
e^{-z\,L}\,\eta_-(z)\,dz,
\end{equation}
where $\Gamma_{\pm}$ is the half ray $\re^+\,e^{\pm
i\,(\pi/2-\theta)}$,
\begin{equation}\label{phi-L:eta}
\eta_{\pm}(z)
=
\frac1{2\,\pi\,i}\,\int_{\gamma_{\pm}}
e^{\zeta\,z}\,\varphi(\zeta)\,d\zeta,
\qquad
z\in\Gamma_{\pm},
\end{equation}
with $\gamma_{\pm}$ being the half-ray $\re^+\,e^{\pm i\,\nu}$ (the
orientation of the paths is not needed in what follows so we do not
pay attention to it). Note that
\begin{equation}
\label{eq:esteta}
|\eta_{\pm}(z)| \lesssim \min ( 1, |z|^{-s-1}),  \qquad
z\in\Gamma_{\pm},
\end{equation}
hence the representation \eqref{phi-L} converges in norm in $\mathcal{L}(L^2)$.
  Usual arguments show the functional property
$\varphi(L)\,\psi(L)=(\varphi\,\psi)(L)$ for two such functions
$\varphi,\psi$.

Any $L$ as above is maximal-accretive and so it has a bounded
holomorphic functional calculus on $L^2$. Given any angle
$\mu\in(\vartheta,\pi)$:
\begin{list}{$(\theenumi)$}{\usecounter{enumi}\leftmargin=.8cm
\labelwidth=0.7cm\itemsep=0.3cm\topsep=.3cm
\renewcommand{\theenumi}{\alph{enumi}}}

\item   For any function $\varphi$, holomorphic and
bounded in $\Sigma_\mu$, the operator $\varphi(L)$ can be defined and is bounded on
$L^2$ with
$$
\|\varphi(L)f\|_{2}
\le
C\,\|\varphi\|_{\infty}\,\|f\|_{2}
$$
where $C$  only depends on $\vartheta$ and $\mu$.

\item For any sequence $\varphi_{k}$ of bounded and holomorphic functions
on $\Sigma_{\mu}$ converging uniformly  on compact subsets of
$\Sigma_{\mu}$ to $\varphi$, we have that $\varphi_{k}(L)$ converges
strongly to $\varphi(L)$ in $\mathcal{L}(L^2)$.

\item The product rule $\varphi(L)\,\psi(L)=(\varphi\,\psi)(L)$ holds for  any two
bounded and holomorphic  functions $\varphi,\psi$  in $\Sigma_{\mu}$.
\end{list}

Let us point out that for more general holomorphic functions (such as
powers), the operators $\varphi(L)$ can be defined as unbounded
operators.

Given a functional Banach space $X$, we say that $L$ has a bounded
holomorphic functional calculus on $X$ if for any
$\mu\in(\vartheta,\pi)$, for any $\varphi$ holomorphic and satisfying \eqref{decay:varphi} in
$\Sigma_\mu$ one has
\begin{equation}
\label{eq:fcX}
\|\varphi(L)f\|_{X}
\le
C\,\|\varphi\|_{\infty}\,\|f\|_{X},
\qquad
f\in X\cap L^2,
\end{equation}
where $C$ depends only on $X$, $\vartheta$ and $\mu$ (but not on the decay of $\varphi$).

If  $X=L^p(w)$ as below,  then \eqref{eq:fcX} implies that
$\varphi(L)$ extends to a bounded operator on $X$ by density. That
$(a)$, $(b)$ and $(c)$ hold with $L^2$ replaced by $X$ for all
bounded holomorphic functions in $\Sigma_{\mu}$, follow from the
theory in \cite{Mc} using the fact that on those $X$, the semigroup
$\{e^{-t\, L}\}_{{t>0}}$ has an infinitesimal generator which is of
type $\vartheta$ (see the last remark of previous section).  We skip
such classical arguments of functional calculi.

\begin{theor}\cite{BK1, Aus}\label{theor:B-K:Aus}
The interior of the set of exponents $p\in(1,\infty)$ such that $L$ has a
bounded holomorphic functional calculus on $L^p$ is equal to
$\Int\J(L)$ defined in Proposition \ref{prop:sgfull}.
\end{theor}

Our first result is  a weighted version of this theorem. We mention
\cite{Mar} where similar weighted estimates are proved under kernel upper
bounds assumptions.

\begin{theor}\label{theor:B-K:weights}
Let $w\in A_\infty$ be such that $\W_w\big(p_-(L),p_+(L)\big)\neq
\emptyset$. Let  $p\in \Int\J_{w}(L)$ and $\mu\in(\vartheta,\pi)$. For any $\varphi$ holomorphic   on $\Sigma_\mu$  satisfying \eqref{decay:varphi}, we have
\begin{equation}
\label{eq:fcw}
\|\varphi(L)f\|_{L^p(w)}
\le
C\,\|\varphi\|_{\infty}\,\|f\|_{L^p(w)},
\qquad
f\in L^\infty_{c},
\end{equation}
with $C$ independent of $\varphi$ and $f$. Hence, $L$ has a bounded holomorphic functional
calculus on $L^p(w)$.
\end{theor}

\begin{remark}\rm Fix $w\in A_{\infty}$  with $\W_w\big(p_-(L),p_+(L)\big)\neq
\emptyset$. If $1<p<\infty$ and  $L$ has a bounded holomorphic functional
calculus on $L^p(w)$, then $p\in \widetilde\J_{w}(L)$. Indeed,  take $\varphi(z)=e^{-z}$. As $\Int\widetilde\J_{w}(L)=\Int\J_{w}(L)$ by Proposition \ref{prop:sgfull}, this shows that range obtained in the theorem is optimal up to endpoints. \end{remark}

\begin{proof} It is enough to assume $\mu<\pi/2$. Note that the operators $e^{-z\, L}$ are uniformly bounded on $L^p(w)$ when $z \in \Sigma_{\mu}$, hence, by \eqref{eq:esteta}, the representation \eqref{phi-L} converges  in norm in $\mathcal{L}(L^p(w))$. Of course, this simple argument does not yield the right estimate, \eqref{eq:fcw}, which is our goal. It is no loss of generality to assume that $\|\varphi\|_{\infty}=1$.
We
split the argument  into three cases: $p\in(\widetilde{p}_-,\widetilde{p}_+)$,
$p\in(\widetilde{p}_-,\widehat{p}_+)$, $p\in(\widehat{p}_-,\widetilde{p}_+)$.

\

\noindent \textit{Case $p\in(\widetilde{p}_-,\widetilde{p}_+)$}.  By $(iii)$ and $(iv)$ in Proposition \ref{prop:weights}, there exist $p_0, q_0$ such that
$$
p_-<p_0<p<q_0<p_+
\qquad {\rm and} \qquad
w\in A_{\frac{p}{p_0}}\cap RH_{\left(\frac{q_0}{p}\right)'}.
$$
The desired bound \eqref{eq:fcw} follows on applying  Theorem
\ref{theor:main-w}  for the underlying measure $dx$ and weight $w$ to
$T=\varphi(L)$ with  $p_{0}, q_{0}$,
$\A_r=I-(I-e^{-r^2\,L})^m$ where $m\ge 1$ is an integer to be chosen
and $S=I$ . As  $\varphi(L)$ and $(I-e^{-r^2\,L})^m$ are bounded on
$L^{p_0}$ (uniformly with respect to $r$ for the latter) by
Proposition \ref{prop:sgfull} and  Theorem \ref{theor:B-K:Aus}, it
remains to checking both \eqref{T:I-A} and \eqref{T:A} on $\D=L^\infty_{c}$.

We start by showing \eqref{T:A}.  We fix  $f\in L^\infty_{c}$ and a
ball $B$. We will use several times the following decomposition of
any given function $h$:
\begin{equation}\label{decomp-h}
h
=
\sum_{j\ge 1} h_j, \qquad\qquad h_{j}=
 h\,\bigchi_{C_j(B)}.
\end{equation}
Fix $1\le k\le m$. Since $p_0\le q_0$ and $p_0, q_0\in  \J(L)$, we have  $e^{-t\,L}\in\off{p_0}{q_0}$ (we are using the equivalence between the two notions of off-diagonal estimates for the Lebesgue measure), hence
\begin{eqnarray*}
\Big( \aver{B}
|e^{-k\,r^2\,L}h_j|^{q_0}\,dx \Big)^{\frac{1}{q_0}}
\lesssim
 2^{j\,(\theta_1+\theta_2)}\,e^{- \alpha\,4^j}\,
\Big(
\aver{2^{j+1}\,B} |h|^{p_0}\,dx \Big)^{\frac{1}{p_0}}
\end{eqnarray*} and by Minkowski's inequality
\begin{eqnarray}\label{eq:T:A}
\Big(\aver{B} |e^{-k\,r^2\,L}h|^{q_0}\,dx
\Big)^{\frac1{q_0}}
\lesssim
\sum_{j\ge 1} g(j)\,
\Big(
\aver{2^{j+1}\,B} |h|^{p_0}\,dx \Big)^{\frac{1}{p_0}}
\end{eqnarray}
with $g(j)=2^{j\,(\theta_1+\theta_2)}\,e^{- \alpha\,4^j}$
for any $h\in L^{p_{0}}$.
This estimate with $h=\varphi(L) f \in L^{p_{0}}$ yields \eqref{T:A} since, by the commutation rule,
$\varphi(L) e^{-k\,r^2\,L}f= e^{-k\,r^2\,L}h$.

We next show \eqref{T:I-A}.  Let $f\in L^\infty_{c}$ and $B$ be
a ball. Write $f=\sum_{j\ge 1} f_j$ as before. For $j=1$, we use the
$L^{p_{0}}$ boundedness of $\varphi(L)$ and $(I-e^{-r^2\,L})^m$,
hence
\begin{equation}\label{varphi-L:w:f1}
\Big(\aver{B} |\varphi(L) (I-e^{-r^2\,L})^m
f_1|^{p_0}\,dx\Big)^{\frac1{p_0}}
\lesssim
\Big(\aver{4\,B} |f|^{p_0}\,dx\Big)^{\frac1{p_0}}.
\end{equation}
For $j\ge 2$, the functions $\eta_{\pm}$
associated with $\psi(z)=\varphi(z)\,(1-e^{-r^2\,z})^m$  by \eqref{phi-L:eta}  satisfy
\begin{equation*}\label{eta:B-K}
|\eta_{\pm}(z)|
\lesssim \frac{r^{2\,m}}{|z|^{m+1}},
\qquad
z\in\Gamma_\pm.
\end{equation*}
Since
$p_0\in \J(L)$,
$\{e^{-z\,L}\}_{z\in\Gamma_{\pm}}\in\off{p_0}{p_0}$ and so
\begin{eqnarray}
\lefteqn{\hskip-1.5cm
\Big( \aver{B}
\Big| \int_{\Gamma_+}\eta_+(z)\, e^{-z\,L} f_j\,dz
\Big|^{p_0}\,dx\Big )^{\frac1{p_0}}
\le
\int_{\Gamma_+} \Big(\aver{B} |e^{-z\,L}
f_j|^{p_0}\,dx\Big)^{\frac1{p_0}}\, |\eta_+(z)|\,|dz|} \nonumber
\\
&\lesssim&
2^{j\,\theta_1} \int_{\Gamma_+}
\dec{\frac{2^j\,r}{\sqrt{|z|}}}^{\theta_2}\,
\expt{-\frac{\alpha\,4^j\,r^2}{|z|}}\,
\frac{r^{2\,m}}{|z|^{m+1}} {|dz|} \,
\Big(\aver{C_j(B)}
|f|^{p_0}\,dx\Big)^{\frac1{p_0}}
\nonumber
\\
&\lesssim&
2^{j\,(\theta_1-2\, m)} \, \Big(\aver{C_j(B)}
|f|^{p_0}\,dx\Big)^{\frac1{p_0}}  \nonumber
\end{eqnarray}
provided $2\, m>\theta_2$. We have used, after a change of variable,
that
$$
\int_0^\infty \dec{s}^{\theta_2}\, e^{-c\,s^2}\,
s^{2\,m}\,\frac{ds}{s} <\infty.
$$
The same is obtained
when one deals with the term corresponding to $\Gamma_-$. Plugging both estimates into
the representation of $\psi(L)$ given by \eqref{phi-L} one obtains
\begin{equation}
\label{eq:ref1}
\Big(\aver{B} |\varphi(L) (I-e^{-r^2\,L})^m
f_j|^{p_0}\,dx\Big)^{\frac{1}{p_0}}
\lesssim
2^{j\,(\theta_1 -2\,m)} \, \Big(\aver{C_j(B)}
|f|^{p_0}\,dx\Big)^{\frac1{p_0}},
\end{equation}
therefore,  \eqref{T:I-A} holds when $2m>\max\{\theta_{1},
\theta_{2}\}$ since $C_{j}(B) \subset 2^{j+1}\, B$.

\

\noindent \textit{Case $p\in(\widetilde{p}_-,\widehat{p}_+)$}: Take
$p_0, q_0$ such that $\widetilde{p}_-<p_0<\widetilde{p}_+$ and
$p_0<p<q_0<\widehat{p}_+$. Let $\A_r=I-(I-e^{-r^2\,L})^m$ for some
large enough $m\ge 1$. Remark that by the previous case, $\varphi(L)$
has the right norm in $\mathcal{L}(L^{p_{0}}(w))$ and so does
$\A_{r}$ by Proposition \ref{prop:sg-w:extension}.  We apply Theorem
\ref{theor:main-w} with the Borel doubling  measure $dw$ and no
weight. Thus, it is enough to see that $\varphi(L)$ satisfies
\eqref{T:I-A} and \eqref{T:A} for $dw$ on $\D=L^\infty_{c}\subset L^{p_{0}}(w)$. But this follows by adapting
the preceding  argument replacing everywhere $dx$ by $dw$ and
observing that $e^{-z\, L} \in \offw{p_{0}}{q_{0}}$ since ${p_{0}},{q_{0}}\in \Int\J_{w}(L)$ and $p_{0}\le q_{0}$. We skip details.

\

\noindent \textit{Case $p\in(\widehat{p}_-,\widetilde{p}_+)$}: Take
$p_0, q_0$ such that $\widetilde{p}_-<q_0<\widetilde{p}_+$ and
$\widehat{p}_-<p_0<p<q_0$. Set $\A_r=I-(I-e^{-r^2\,L})^m$ for some
integer $m\ge 1$ to be chosen later. Since $q_0\in (\widetilde
p_-,\widetilde p_+)$,  by the first case, $\varphi(L)$ has the right
norm in $\mathcal{L}(L^{q_{0}}(w))$ and so does $\A_{r}$ by
Proposition \ref{prop:sg-w:extension}.  We  apply Theorem
\ref{theor:SHT:small} with underlying measure $dw$.  It is enough to
show \eqref{SHT:small:T:I-A} and \eqref{SHT:small:A}.  Fix a ball $B$
and $f \in L^\infty_{c}$  supported in $B$.

We begin with \eqref{SHT:small:A} for $\A_{r}$. It is
enough to show it for $e^{-k\,r^2\,L}$ with $1\le k\le m$. Since $p_{0}, q_{0}\in \J_{w}(L)$ and $p_{0}\le q_{0}$ we have
$e^{-t\,L}\in\offw{p_0}{q_0}$, hence
\begin{equation}\label{est:Ar:w}
\Big(
\aver{C_j(B)} |e^{-k\,r^2\,L} f|^{q_0}\,dw
\Big)^{\frac1{q_0}}
\lesssim
2^{j\,(\theta_1+\theta_2)}\,e^{-c\,4^j} \Big(\aver{B}
|f|^{p_0}\,dw\Big)^{\frac1{p_0}}.
\end{equation}
This implies \eqref{SHT:small:A} with
$g(j)=C\,2^{j\,(\theta_1+\theta_2)}\,e^{-c\,4^j}$ and
$\sum_{j\ge 1}g(j)\,2^{D\,j}<\infty$ holds where $D$ is the doubling order of $dw$.

We turn to  \eqref{SHT:small:T:I-A}.   Let $j\ge 2$.   The argument is the same as the one for \eqref{eq:ref1} by reversing the roles of $C_{j}(B)$ and $B$, and using $dw$ and $e^{-z\,L}\in\offw{p_0}{p_0}$ (since $p_{0}\in \J_{w}(L)$)  instead of $dx$ and $e^{-z\,L}\in\off{p_0}{p_0}$. We obtain
\begin{equation*}
\label{eq:ref1bis}
\Big(\aver{C_{j}(B)} |\varphi(L) (I-e^{-r^2\,L})^m
f|^{p_0}\,dw\Big)^{\frac{1}{p_0}}
\lesssim
2^{j\,(\theta_1 -2\,m)} \, \Big(\aver{B}
|f|^{p_0}\,dw\Big)^{\frac1{p_0}}
\end{equation*}
provided $2\, m >\theta_{2}$ and it remains to impose further $2\, m > \theta_{1}+D$  to conclude.
\end{proof}

\begin{remark}
\rm If $\W_{w}(p_{-}, p_{+})\ne \emptyset$, the last part of the
proof yields  weighted weak-type $ (\widehat{p}_{-},
\widehat{p}_{-})$ of $\varphi(L)$  provided $\widehat{p}_{-}\in
\J_{w}(L)$. To do so on only has to take $p_{0}= \widehat{p}_{-}$.
\end{remark}

\section{Riesz transforms}\label{sec:RT}

The Riesz transforms associated to $L$ are $\partial_{j}L^{-1/2}$, $1\le j \le n$. Set   $\nabla L^{-1/2}=(\partial_{1}L^{-1/2}, \ldots, \partial_{n}L^{-1/2})$. The
solution of the Kato conjecture \cite{AHLMcT} implies that this operator extends boundedly to $L^2$ (we ignore the $\co^n$-valued aspect of things). This allows the representation
\begin{equation}
\label{eq:RT}
\nabla L^{-1/2}f = \frac{1}{\sqrt \pi}\int_{0}^\infty \nabla e^{-t\, L}f\, \frac{dt}{\sqrt t},
\end{equation}
in which the integral converges  strongly in $L^2$ both at $0$ and
$\infty$ when $f\in L^2$. Note that for an arbitrary $f\in L^2$,
$h=L^{-1/2}f$ makes sense in the homogeneous Sobolev space $\dot
H^{1}$ which is the completion of $C_{0}^\infty(\re^n)$ for the
semi-norm $\|\nabla h\|_{2}$ and $\nabla$ becomes the extension  of
the gradient to that space. This construction can be forgotten if
$n\ge 3$ as $\dot H^{1}\subset L^{2^*}$ but not if $n\le 2$. To
circumvent this technical difficulty, we  introduce $S_{\ep}=
\frac{1}{\sqrt \pi}\int_{\ep}^{1/\ep}  e^{-t\, L}\, \frac{dt}{\sqrt
t}$ for $0<\ep<1$ and, in fact, $\nabla S_{\ep}$ are uniformly
bounded on $L^2$ and converge strongly in $L^2$. This defines $\nabla
L^{-1/2}$.

\begin{theor}[\cite{Aus}]\label{theor:Riesz-Auscher}
The maximal interval  of exponents $p\in(1,\infty)$ for which $\nabla L^{-1/2}$ has a  bounded extension to  $L^p$ is equal to  $\Int\K(L)$
defined in Proposition \ref{prop:sgfull} and for $p\in \Int\K(L)$,
$\|\nabla f \|_{p}\sim \|L^{1/2}f\|_{p}$ for all $f\in
\D(L^{1/2})=H^{1}$ \textup{(}the Sobolev space\textup{)}.
\end{theor}

Again, the operators $\nabla S_{\ep}$ are uniformly bounded on $L^p$ and converge strongly in $L^p$ as $\ep \to 0$. Indeed, for $f\in L^\infty_{c}$, $S_{\ep}f\in \D(L)\subset \D(L^{1/2})$ and
$
\|\nabla S_{\ep}f\|_{p} \lesssim \|L^{1/2}S_{\ep}f\|_{p}.
$
Observe that $L^{1/2}S_{\ep}=\varphi_{\ep}(L)$, where $\varphi_{\ep}$
is a bounded holomorphic function in $\Sigma_{\mu}$ for any
$0<\mu<\pi/2$ with $\sup _{\ep}\|\varphi_{\ep}\|_{\infty} <\infty$
and $\{\varphi_{\ep}\}_\ep$ converges uniformly to 1 on compact
subsets of $\Sigma_{\mu}$ as $\ep\to 0$. The claim follows by Theorem
\ref{theor:B-K:Aus} and density.

We turn to  weighted norm inequalities. Remark that by Proposition \ref{prop:sg-w:extension}, for all
$p \in \J_w(L)$, $S_{\ep}$ is bounded on $L^p(w)$ (the norm must depend on $\ep$) and for all $p\in
\K_{w}(L)$, $\nabla S_{\ep}$ is bounded on $L^p(w)$ with no control yet on  the norm with respect to $\ep$.

\begin{theor}\label{theor:ext-RT}
Let $w\in A_\infty$ be such that $\W_w\big(q_-(L),q_+(L)\big)\neq
\emptyset$. For all $p\in
\Int \K_w(L)$ and $ f \in L^\infty_{c}$,
\begin{equation}
\label{eq:Riesz-w } \|\nabla L^{-1/2} f\|_{L^p(w)} \lesssim
\|f\|_{L^p(w)}.
\end{equation}
Hence,  $\nabla L^{-1/2}$ has a bounded extension to  $L^p(w)$.
\end{theor}

We note that  for a given $p \in \Int\K_{w}(L)$, once \eqref{eq:Riesz-w } is established, similar arguments  using
Theorem  \ref{theor:B-K:weights} imply
 convergence
in $L^p(w)$ of $\nabla S_{\ep}f$ to $\nabla L^{-1/2}f$ for $f\in L^\infty_{c}$.

\begin{proof}

We split the argument in
three  cases: $p\in(\widetilde{q}_-,\widetilde{q}_+)$, $p\in(\widetilde{q}_-,\widehat{q}_+)$,
$p\in(\widehat{q}_-,\widetilde{q}_+)$.

\

\noindent \textit{Case $p\in(\widetilde{q}_-,\widetilde{q}_+)$}:   By $(iii)$ and $(iv)$ in
Proposition \ref{prop:weights}, there exist $p_0, q_0$ such that
\begin{equation*}\label{exponents}
q_-<p_0<p<q_0<q_+
\qquad {\rm and} \qquad
w\in A_{\frac{p}{p_0}} \cap RH_{\left(\frac{q_0}{p}\right)'}.
\end{equation*}
The desired estimate \eqref{eq:Riesz-w } is obtained by applying
Theorem \ref{theor:main-w} with  underlying measure $dx$ and weight
$w$. Hence, it suffices  to verify \eqref{T:I-A} and \eqref{T:A}  on
$\D=L^\infty_{c}$ for $T=\nabla L^{-1/2}$, $S=I$ and
$\A_r=I-(I-e^{-r^2\,L})^m$, with $m$ a  large enough integer.  These
conditions will be proved as in \cite{Aus}, but here we use the whole
range of exponents for which the Riesz transforms are bounded on
unweighted $L^p$ spaces, that is, $(q_{-}, q_{+})$.

 \begin{lemma}\label{lemma:est-Riesz}
Fix a  ball $B$. For $f\in L^\infty_{c}$ and $m$ large enough,
\begin{equation}
\label{eq:est-Riesz:1}
\Big( \aver{B} |\nabla
L^{-1/2}(I-e^{-r^2\,L})^mf|^{p_0}\, dx
\Big)^{\frac1{p_0}}
\le
\sum_{j\ge 1} g_1(j)
\Big( \aver{C_{j}(B)}\!\!
|f|^{p_0}\, dx\Big)^{\frac1{p_0}}
\end{equation}
and for $f \in L^{p_{0}}$ such that $\nabla f \in L^{p_{0}}$ and  $1\le k\le m$,
\begin{equation}
\label{eq:est-Riesz:2}
\Big( \aver{B} |\nabla
e^{-k\,r^2\,L}f|^{q_0}\,dx\Big)^{\frac1{q_0}}
\le
\sum_{j\ge 1} g_2(j)\,
\Big( \aver{2^{j+1}\,B}
|\nabla f|^{p_0}\,dx\Big)^{\frac1{p_0}},
\end{equation}
where $g_1(j)=C_{m}\,2^{j\,\theta}\,4^{-m\,j}$ and
$g_2(j)=C_m\,2^j\,\sum_{l\ge j} 2^{l\,\theta}\,e^{-\alpha\,4^l}$ for some $\theta>0$.
\end{lemma}

Assume this is proved. Note that  if $2m>\theta$ then $\sum_{j\ge
1}g_1(j)<\infty$ and  the first estimate is \eqref{T:I-A}.

Next, expanding $\A_r=I-(I-e^{-r^2\,L})^m$, the latter estimate
applied to  $S_{\ep}f$  in place of $f$ (since  $S_{\ep}f\in
L^{p_{0}}$ and $\nabla S_{\ep}f \in L^{p_{0}}$)   and the commuting
rule $\A_{r}S_{\ep}=S_{\ep} \A_{r}$ give us
$$
\Big( \aver{B} |\nabla S_{\ep} \A_{r}
f|^{q_0}\,dx\Big)^{\frac1{q_0}}
\le
\sum_{j\ge 1} g_2(j)\,
\Big( \aver{2^{j+1}\,B}
|\nabla S_{\ep} f|^{p_0}\,dx\Big)^{\frac1{p_0}}.
$$
By letting $\ep$ go to 0 (the justification uses the observations made at the beginning of this section and is left to the reader), we obtain
 \eqref{T:A} using $\sum_{j\ge 1}g_2(j)<\infty$.
Therefore, by  Theorem \ref{theor:main-w}, \eqref{eq:Riesz-w } holds
for $f\in L^\infty_{c} $.

\begin{proof}[Proof of Lemma \ref{lemma:est-Riesz}]
 We begin with the first estimate. Decomposing $f$ as in
\eqref{decomp-h},
\begin{eqnarray*}
\Big(\aver{B} |\nabla L^{-1/2} (I-e^{-r^2\,L})^m
f|^{p_0}\,dx\Big)^{\frac1{p_0}}
\le
\sum_{j\ge 1}
\Big(
\aver{B} |\nabla L^{-1/2} (I-e^{-r^2\,L})^m f_j|^{p_0}\,dx
\Big)^{\frac1{p_0}}.
\end{eqnarray*}
For $j=1$, since $q_-<p_0<q_+$,  $\nabla
L^{-1/2}$ and $e^{-r^2\,L}$ are bounded on $L^{p_0}$ by Theorem \ref{theor:Riesz-Auscher}
and Proposition \ref{prop:sgfull}. Hence,
$$
\Big( \aver{B} |\nabla L^{-1/2}(I-e^{-r^2\,L})^m
f_1|^{p_0}\,dx \Big)^{\frac1{p_0}}
\lesssim
\Big( \aver{4\,B} |f|^{p_0}\,dx\Big)^{\frac1{p_0}}.
$$
For $j\ge 2$, we use a different approach. If $h\in L^2$, by
\eqref{eq:RT}
\begin{equation*}\label{RT-repre}
\nabla L^{-1/2}(I-e^{-r^2\,L})^m h
=
\frac1{\sqrt{\pi}} \int_0^\infty \sqrt{t}\,\nabla \varphi(L,t)
h\,\frac{dt}{t},
\end{equation*}
where $\varphi(z,t)=e^{-t\,z}\,(1-e^{-r^2\,z})^m$. The functions
$\eta_{\pm}(\cdot, t)$ associated with $\varphi(\cdot, t)$ by
\eqref{phi-L:eta}  satisfy
\begin{equation*}\label{eta-RT}
|\eta_{\pm}(z,t)|
\lesssim
\frac{r^{2\,m}}{(|z|+t)^{m+1}},
\qquad
z\in \Gamma_{\pm}, t>0.
\end{equation*}
Since $\sqrt{z}\,\nabla e^{-z\,L}\in\off{p_0}{p_0}$ (note
that $p_0\in  \K(L)$ and we are using the equivalence between the two notions of off-diagonal estimates for the Lebesgue measure),
\begin{eqnarray*}
\lefteqn{
\Big( \aver{B}
\Big|
\int_{\Gamma_+} \eta_+(z)\, \sqrt{t}\,\nabla e^{-z\,L} f_j\,dz
\Big|^{p_0}\,dx\Big)^{\frac1{p_0}}} \nonumber
\\
&\le&
\int_{\Gamma_+}
\Big(
\aver{B} |\sqrt{z}\,\nabla e^{-z\,L}
f_j|^{p_0}\,dx\Big)^{\frac1{p_0}}\,
\frac{\sqrt{t}}{\sqrt{|z|}}\,|\eta_+(z)|\,|dz| \nonumber
\\
&\lesssim&
2^{j\,\theta_1} \int_{\Gamma_+}
\dec{\frac{2^j\,r}{\sqrt{|z|}}}^{\theta_2}
\expt{-\frac{\alpha\,4^j\,r^2}{|z|}}\,
\frac{\sqrt{t}}{\sqrt{|z|}}\,|\eta_+(z)|\,|dz| \,
\Big(\aver{C_j(B)}
|f|^{p_0}\,dx\Big)^{\frac1{p_0}} \nonumber
\\
&\lesssim&
2^{j\,\theta_1} \int_0^\infty
\dec{\frac{2^j\,r}{\sqrt{s}}}^{\theta_2}
\expt{-\frac{\alpha\,4^j\,r^2}{s}}\, \frac{\sqrt{t}}{\sqrt{s}}
\,\frac{r^{2\,m}}{(s+t)^{m+1}}\,ds\, \Big(\aver{C_{j}(B)}
|f|^{p_0}\,dx\Big)^{\frac1{p_0}}.
\end{eqnarray*}
Observing that when  $2\, m>\theta_{2}$
\begin{eqnarray*}
\int_{0}^\infty\int_0^\infty
\dec{\frac{2^j\,r}{\sqrt{s}}}^{\theta_2}
\expt{-\frac{\alpha\,4^j\,r^2}{s}}\, \frac{\sqrt{t}}{\sqrt{s}}\,\frac{r^{2\,m}}{(s+t)^{m+1}}\,ds\frac {dt}t = C\, 4^{-j\, m},
\end{eqnarray*}
and plugging this, plus the corresponding term for $\Gamma_-$,  into
the representation \eqref{phi-L}, we obtain
\begin{eqnarray}
\Big( \aver{B} |\nabla e^{-t\,L}(I-e^{-r^2\,L})^m
  f_j
|^{p_0}\,dx\Big)^{\frac1{p_0}}
&\lesssim&
\int_0^\infty
\Big( \aver{B}
|\sqrt{t}\,\nabla \varphi(L,t) f_j|^{p_0}\,dx\Big)^{\frac1{p_0}} \, \frac{dt}{t} \nonumber
\\
&\lesssim& 2^{j\,(\theta_1-2\,m)}
\Big(\aver{C_{j}(B)}
|f|^{p_0}\,dx\Big)^{\frac1{p_0}}. \label{estimate:RT-ext}
\end{eqnarray}
This readily yields the first estimate in the lemma.

Let us get the second one. Fix $1\le k\le m$. Let $f\in L^{p_{0}}$ such that $\nabla f \in L^{p_{0}}$. We write
$h=f-f_{4\,B}$ where $f_{\lambda\,B}$ is the $dx$-average of $f$ on $\lambda
\,B$. Then by the conservation property (see \cite{Aus})
$e^{-t\, L}1=1$ for all
 $t>0$,
 we have
$$
\nabla e^{-k\,r^2\,L} f
=
\nabla e^{-k\,r^2\,L} (f-f_{4\,B})
=
\nabla e^{-k\,r^2\,L} h
=
\sum_{j\ge 1} \nabla e^{-k\,r^2\,L}h_j,
$$
with $h_j=h\,\bigchi_{C_j(B)}$. Hence,
$$
\Big(\aver{B} |\nabla
e^{-k\,r^2\,L}f|^{q_0}\,dx\Big)^{\frac1{q_0}}
\le
\sum_{j\ge 1}
\Big(\aver{B}  |\nabla
e^{-k\,r^2\,L}h_j|^{q_0}\,dx\Big)^{\frac1{q_0}}.
$$
Since $p_0\le  q_0$ and $p_0, q_0\in  \K(L)$, $\sqrt{t}\,\nabla
e^{-t\,L}\in\off{p_0}{q_0}$. This and the $L^{p_{0}}$-Poincar\'e
inequality for $dx$ yield
\begin{eqnarray}\label{eq:poincare}
\lefteqn{\hskip-0.5cm
\Big(\aver{B}
|\nabla e^{-k\,r^2\,L}h_j|^{q_0}\,dx\Big)^{\frac1{q_0}}
\lesssim
\frac{2^{j\,(\theta_1+\theta_2)}\,e^{-\alpha\,4^j}}r\,
\Big(\aver{C_j(B)}
|h_j|^{p_0}\,dx\Big)^{\frac1{p_0}}}
\nonumber
\\
&\le&
\frac{2^{j\,(\theta_1+\theta_2)}\,e^{-\alpha\,4^j}}r\,
\Big(\aver{2^{j+1}\,B}|f-f_{4\,B}|^{p_0}\,dx\Big)^{\frac1{p_0}}
\nonumber
\\
&\le&
\frac{2^{j\,(\theta_1+\theta_2)}\,e^{-\alpha\,4^j}}r\,
\bigg(
\Big(\aver{2^{j+1}\,B}\hskip-7pt
|f-f_{2^{j+1}\,B}|^{p_0}\,dx\Big)^{\frac1{p_0}} +\sum_{l=2}^j
|f_{2^l\,B}-f_{2^{l+1}\,B}| \bigg)
\nonumber
\\
&\lesssim&
\frac{2^{j\,(\theta_1+\theta_2)}\,e^{-\alpha\,4^j}}r\,\sum_{l=1}^{j}
\Big(\aver{2^{l+1}\,B}
|f-f_{2^{l+1}\,B}|^{p_0}\,dx\Big)^{\frac1{p_0}}
\nonumber
\\
&\lesssim&
2^{j\,(\theta_1+\theta_2)}\,e^{-\alpha\,4^j}\,\sum_{l=1}^{j} 2^{l}\,
\Big(\aver{2^{l+1}\,B}
|\nabla f|^{p_0}\,dx\Big)^{\frac1{p_0}},
\end{eqnarray}
which is the desired estimate with $\theta=\theta_{1}+\theta_{2}$.
\end{proof}

\noindent \textit{Case $p\in(\widetilde{q}_-,\widehat{q}_+)$}: Take
$p_0, q_0$ such that $\widetilde{q}_-<p_0<\widetilde{q}_+$ and
$p_0<p<q_0<\widehat{q}_+$. Let $\A_r=I-(I-e^{-r^2\,L})^m$ for some
$m\ge 1$ to be chosen later. As   $p_0\in (\widetilde q_-, \widetilde
q_+)$,  both $\nabla L^{-1/2}$ and $\A_{r}$  are bounded on
$L^{p_0}(w)$ (we have just shown it for the Riesz transforms and
$\A_r$ are bounded uniformly in $r$ by Proposition
\ref{prop:sg-w:extension}). By Theorem \ref{theor:main-w} with
underlying doubling measure $dw$ and no weight, it is enough to
verify \eqref{T:I-A} and \eqref{T:A} on $\D=L^\infty_{c}$ for
$T=\nabla L^{-1/2}$, $S=I$ and $\A_{r}$.  To do so, it suffices to
copy the proof of Lemma \ref{lemma:est-Riesz} in the weighted case by
changing systematically $dx$ to $dw$, off-diagonal estimates with
respect to $dx$ by those with respect to $dw$ given the choice of
$p_{0},q_{0}$. Also in the argument with $dx$ we used a Poincar\'e
inequality. Here, since $p_0\in \W_{w}(q_-,q_+)$,  $w\in A_{p_0/q_-}$
and  in particular $w\in A_{p_0}$ (since $q_-\ge 1$). Therefore we
can use the $L^{p_0}(w)$-Poincar\'e inequality (see \cite{FPW}):
$$
\Big( \aver{B} |f-f_{B,w}|^{p_{0}}\, dw\Big)^{\frac 1 {p_{0}}} \lesssim r(B) \Big( \aver{B} |\nabla f|^{p_{0}}\, dw\Big)^{\frac 1 {p_{0}}}
$$
for all $f\in L_{\rm loc}^1(w)$ such that $\nabla f \in
L^{p_{0}}_{\rm loc}(w)$, where $f_{B,w}$ is the $dw$-average of $f$
over $B$. We leave further details to the reader.

\

\noindent \textit{Case $p\in(\widehat{q}_-,\widetilde{q}_+)$}: Take
$p_0, q_0$ such that $\widetilde{q}_-<q_0<\widetilde{q}_+$ and
$\widehat{q}_-<p_0<p<q_0$. Set $\A_r=I-(I-e^{-r^2\,L})^m$ for some
integer $m\ge 1$ to be chosen later. Since $q_0\in (\widetilde q_-,
\widetilde q_+)$, it follows that $\nabla L^{-1/2}$ is already
bounded on $L^{q_0}(w)$ and so is $\A_r$.  That $\nabla L^{-1/2}$ is
 bounded on $L^{p}(w)$ will follow on applying Theorem \ref{theor:SHT:small} with underlying measure $w$. Hence it is enough to check both \eqref{SHT:small:T:I-A} and  \eqref{SHT:small:A}.

 We begin with
  \eqref{SHT:small:A}. By Proposition
\ref{prop:sg-w:extension}, $\inf \J_w(L)=\inf \K_w(L)=\widehat{q}_-$.
Since $p_0>\widehat{q}_-$ and $p_0 \le q_0\in \W_w(q_-,q_+)\subset
\W_w(p_-,p_+)\subset \J_w(L)$, we have $p_0, q_0\in \J_w(L)$ and so
$e^{-t\,L}\in\offw{p_0}{q_0}$. This yields \eqref{est:Ar:w}, hence
\eqref{SHT:small:A} with
$g(j)=C\,2^{j\,(\theta_1+\theta_2)}\,e^{-c\,4^j}$ which clearly
satisfies $\sum_j g(j)\,2^{D\,j}<\infty$, with $D$ the doubling order
of $dw$.

We next show \eqref{SHT:small:T:I-A}.  Let $f\in L^\infty_{c}$ be
supported on a ball $B$ and $j\ge 2 $. The argument is the same as
the one for \eqref{estimate:RT-ext} by reversing the roles of
$C_{j}(B)$ and $B$, and using $dw$ and $\sqrt z \, \nabla
e^{-z\,L}\in\offw{p_0}{p_0}$ (since $p_{0}\in \K_{w}(L)$)  instead of
$dx$ and $\sqrt z \, \nabla e^{-z\,L}\in\off{p_0}{p_0}$. Hence, we
obtain
\begin{equation*}
\label{eq:ref2bis}
\Big(\aver{C_{j}(B)} |\nabla L^{-1/2} (I-e^{-r^2\,L})^m
f|^{p_0}\,dw\Big)^{\frac{1}{p_0}}
\lesssim
2^{j\,(\theta_1 -2\,m)} \, \Big(\aver{B}
|f|^{p_0}\,dw\Big)^{\frac1{p_0}}
\end{equation*}
provided $2\, m >\theta_{2}$ and it remains to impose further $2\, m > \theta_{1}+D$  to conclude.
\end{proof}

\begin{remark}
\rm If $\W_{w}(q_{-}, q_{+})\ne \emptyset$, the last part of the
proof yields  weighted weak-type $ (\widehat{q}_{-},
\widehat{q}_{-})$ of $\nabla L^{-1/2}$  provided $\widehat{q}_{-}\in
\K_{w}(L)$, one only needs to take  $p_{0}= \widehat{q}_{-}$.
\end{remark}

\begin{remark}
\rm
Theorem \ref{theor:Riesz-Auscher} asserts that  $\Int\K(L)$ is the \emph{exact} range of $L^p$ boundedness for the Riesz transforms when $w=1$. When $w\ne 1$,  we cannot repeat the same argument as it used Sobolev embedding which has no simple counterpart in the weighted situation. However, if we insert in the integral of \eqref{eq:RT} a function $m(t)$ with $m\in L^\infty(0,\infty)$, then
\eqref{eq:Riesz-w } holds with a constant proportional to $\|m\|_{\infty}$.  Indeed,
Let $\varphi_{m}(z)=\int_{0}^\infty z^{1/2} e^{-t\, z}\, m(t)\, \frac{dt}{\sqrt t}$ for $z\in \Sigma_{\mu}$, $\vartheta<\mu<\pi/2$. Then, $\varphi_{m}$ is holomorphic in $\Sigma_{\mu}$ and bounded with $\|\varphi_{m}\|_{\infty}\le c_{\mu}\|m\|_{\infty}$. Now for $f\in L^\infty_{c}$,
$$
\int_{0}^\infty \nabla e^{-t\, L}f\, m(t)\,  \frac{dt}{\sqrt t} = \nabla L^{-1/2} \varphi_{m}(L)f
$$
hence, combining Theorems \ref{theor:B-K:weights} and \ref{theor:ext-RT}, we obtain for $p\in \Int\K_{w}(L)$,
$$
\Big\| \int_{0}^\infty \sqrt t\, \nabla e^{-t\, L}f\, m(t)\,  \frac{dt}{t} \Big\|_{L^p(w)} \lesssim \|m\|_{\infty}\, \|f\|_{L^p(w)}.
$$
Conversely, given an exponent $p\in (1,\infty)$, assume that this
$L^p(w)$ estimate holds for all $m\in L^\infty$. Using randomization
techniques which we skip (see Section \ref{sec:vv} for some account
on such techniques),  this implies
$$
\Big\|\Big(
\int_0^\infty |\sqrt t\,  \nabla e^{-t\,L}f|^2\,\frac{dt}t
\Big)^{\frac12}
\Big\|_{L^p(w)} \lesssim \|f\|_{L^p(w)}.
$$
This square function estimate is proved directly in Section \ref{sec:sf} and we indicate at the end  of that section why  this  inequality implies $p\in \overline{\widetilde\K_{w}(L)}$.  Thus, the range in $p$ is sharp up to endpoints (see Proposition \ref{prop:sg-w:extension}).
\end{remark}

\section{Reverse inequalities for square roots}\label{sec:RRT}

We continue on square roots by studying when the inequality opposite
to \eqref{eq:Riesz-w }  hold. First we recall the unweighted case.

\begin{theor}\label{theor:reverseRiesz}\cite{Aus}
If $\max\big\{1,\frac{n\,p_{-}(L)}{n+p_{-}(L)}\big\}<p<  p_+(L)$ then for $f\in \cals$,
\begin{equation}
\label{eq:reverseRiesz}
\|L^{1/2}f\|_{p} \lesssim \|\nabla f\|_{p}.
\end{equation}
\end{theor}

To state our result, we need  a new exponent.  For $p>0$, define
$$
p_{w,*}
=
\frac{n\, r_{w}\, p}{\, n\, r_{w}+ p},
$$
where $r_{w}=\inf\{r\ge 1\ : \ w \in A_{r} \}$. Set also $p_{w}^*=
\frac{n\, r_{w}\, p}{\, n\, r_{w}- p}$ for $p<n\,r_w$ and
$p_*=\infty$ otherwise. Note that $(p_{w,*})_{w}^*=p$.

\begin{theor}\label{theor:reverseRiesz-w-2}
Let  $w\in A_\infty$ with $\W_w\big(\,p_{-}(L),p_{+}(L)\big)\ne
\emptyset$. If $\max\big\{\,r_{w}\,,\, (\widehat p_{-})_{w,*}\,\big\}<p< \widehat p_+$ then for $f\in \cals$,
\begin{equation}
\label{eq:reverseRiesz:w}
\|L^{1/2}f\|_{L^p(w)} \lesssim \|\nabla f\|_{L^p(w)}.
\end{equation}
\end{theor}

\begin{remark}\label{remark:reverse:dim1-2} \rm
Recall that $ \widetilde p_{-} = p_{-}(L)\,r_{w}$ and we have
$(\widehat p_{-})_{w,*}< \widehat p_{-}\le \widetilde p_{-}$. If
$p_{-}(L)=1$, then $(\widehat p_{-})_{w,*}\le r_w$, so $\max\big\{\,
r_{w}\,,\, (\widehat p_{-})_{w,*}\, \big\}=r_w= \widetilde p_{-}$.
This happens for example  when $L$ is real or when $n=1,2$.
\end{remark}

Define $\dot W^{1,p}(w)$ as the completion of $\cals$ under
the semi-norm $\|\nabla f\|_{L^p(w)}$. Arguing as in \cite{AT1} (see
\cite{Aus}) combining Theorems    \ref{theor:ext-RT}  and
\ref{theor:reverseRiesz-w-2}, we obtain the following consequence.

\begin{corol} Assume $ \W_w\big(\,q_{-}(L),q_{+}(L)\big)\ne \emptyset$.  If $p\in \Int\K_{w}(L)$ with $p>r_{w}$, then $L^{1/2}$
extends to an isomorphism from $\dot W^{1,p}(w)$ into $L^p(w)$.
\end{corol}

\begin{proof}[Proof of Theorem \ref{theor:reverseRiesz-w-2}]
We split the argument in
three  cases: $p\in(\widetilde{p}_-,\widetilde{p}_+)$, $p\in(\widetilde{p}_-,\widehat{p}_+)$,
$p\in(\max\big\{r_{w}\, ,\, (\widehat p_{-})_{w,*}\big\},\widetilde{p}_+)$.

\

\noindent \textit{Case $p\in(\widetilde{p}_-,\widetilde{p}_+)$}:
It relies on the following lemma.

\begin{lemma}\label{lemma:est-reverseRiesz} Let $p_{0}\in \Int\J(L)$ and $q_{0}\in \J(L)$
with $p_{0}<q_{0}$. Let $B$ be a ball and $m\ge 1$ an integer. For all $f\in \cals$, we have
\begin{equation}\label{Tep-1}
\Big( \aver{B} |L^{1/2}(I-e^{-r^2\,L})^mf|^{p_0}\, dx\Big)^{\frac1{p_0}}
\le \sum_{j\ge 1} g_1(j)\,
\Big( \aver{2^{j+1}\,B}
|\nabla f|^{p_0}\,dx\Big)^{\frac1{p_0}}
\end{equation}
 for $m$ large enough depending on
$p_{0}$ and $q_{0}$, and
\begin{equation}\label{Tep-2}
\Big( \aver{B} |L^{1/2} (I-(I-e^{-r^2\,L})^m)f|^{q_0}\,dx\Big)^{\frac1{q_0}}
\le\sum_{j\ge 1} g_2(j)\,
\Big( \aver{2^{j+1}\,B}
|L^{1/2} f|^{p_0}\,dx\Big)^{\frac1{p_0}},
\end{equation}
where $g_1(j)=C_{m}\,2^{j\,\theta}\,4^{-m\,j}$ and
$g_2(j)=C_{m}\,2^{j\,\theta}\,e^{-\alpha\,4^j}$ for some $\theta>0$, and the
implicit constants are independent of $B$ and $f$.
\end{lemma}

Admit this lemma for a moment. Since
$p\in(\widetilde{p}_-,\widetilde{p}_+)=\W_w\big(p_-, p_+\big)$,
by $(iii)$ and $(iv)$ in Proposition \ref{prop:weights}, there exist
$p_0, q_0$ such that
$$
p_-<p_0<p<q_0<p_+
\qquad {\rm and}\qquad
w\in A_{\frac{p}{p_0}}\cap RH_{\left(\frac{q_0}{p}\right)'}.
$$
Note that \eqref{Tep-1} and \eqref{Tep-2} are respectively the
conditions \eqref{T:I-A} and \eqref{T:A} of  Theorem
\ref{theor:main-w} with underlying measure $dx$ and weight $w$, $T=
L^{1/2}$, $\A_r=I-(I-e^{-r^2\,L})^m$, with $m$ large enough, and $S
f=\nabla f$. Hence  we obtain \eqref{eq:reverseRiesz:w}.

\begin{proof}[Proof of Lemma \ref{lemma:est-reverseRiesz}.]
We first show \eqref{Tep-2}. Using the commutation rule and expanding
$(I-e^{-r^2\, L})^m$ it suffices to apply \eqref{eq:T:A}  as $p_{0},q_{0}\in \J(L)$ to
$h=L^{1/2}f$.

We turn to \eqref{Tep-1}.  If $\varphi(z)=z^{1/2} (1-e^{-r^2\, z})^m$, then $\varphi(L)f=L^{1/2}(I-e^{-r^2\,L})^mf$. By
the conservation property
\begin{eqnarray}\label{eq:rep}
\varphi(L)\, f=
\varphi(L)\, (f-f_{4\,B})
=
\sum_{j\ge 1} \varphi(L)\, h_j,
\end{eqnarray}
where $h_{j}=(f-f_{4\,B})\,\phi_j$. Here,
$\phi_j=\bigchi_{C_j(B)}$ for $j\ge 3$, $\phi_1$ is a  smooth
function with support in $4\,B$, $0\le \phi_1\le 1$, $\phi_{1}=1$ in $2\,B$
and $\|\nabla \phi_1\|_{\infty}\le C/r$ and, eventually, $\phi_2$ is
taken so that $\sum_{j\ge 1} \phi_j=1$. We estimate each term in
turn. For $j=1$, since $ p_-<p_0< p_+$, by the bounded holomorphic
functional calculus on $L^{p_{0}}$ (Theorem \ref{theor:B-K:Aus}) and
$\varphi(L)\, h_{1}= (I-e^{-r^2\, L})^m\, L^{1/2}h_{1}$,   one
has uniformly in $r$,
$$
\left\|\varphi(L)\, h_{1}\right\|_{{p_{0}}}
 \lesssim
\|L^{1/2}h_{1}\|_{{p_{0}}}.
$$
Next,  Theorem \ref{theor:reverseRiesz}, $L^{p_{0}}$-Poincar\'e
inequality and the definition of $h_{1}$ imply
$$
\|L^{1/2}h_{1}\|_{{p_{0}}} \lesssim \|\nabla
h_{1}\|_{{p_{0}}}
\lesssim
\|\nabla f\|_{L^{p_{0}}(4\, B)}.
$$
Therefore,
$$
\Big( \aver{B} |\varphi(L)\, h_{1}|^{p_0}\,dx\Big)^{\frac1{p_0}}
\lesssim
\Big( \aver{4\,B} |\nabla f
|^{p_0}\,dx\Big)^{\frac1{p_0}}.
$$
For $j\ge 3$, the functions  $\eta_{\pm}$
associated with $\varphi$  by \eqref{phi-L:eta}  satisfy
$$
|\eta_{\pm}(z)|
\lesssim \frac{r^{2\,m}}{|z|^{m+3/2}},
\qquad
z\in\Gamma_\pm.
$$
Since $p_0\in \J(L)$,
$\{e^{-z\,L}\}_{z\in\Gamma_{\pm}}\in\off{p_0}{p_0}$ and so
\begin{eqnarray*}
\lefteqn{\hskip-1cm
\Big( \aver{B}
\Big| \int_{\Gamma_+}\eta_+(z)\, e^{-z\,L} h_j\,dz
\Big|^{p_0}\,dx\Big )^{\frac1{p_0}}
\le
\int_{\Gamma_+} \Big(\aver{B} |e^{-z\,L}
h_j|^{p_0}\,dx\Big)^{\frac1{p_0}}\, |\eta_+(z)|\,|dz|}
\\
&\lesssim&
2^{j\,\theta_1} \int_{\Gamma_+}
\dec{\frac{2^j\,r}{\sqrt{|z|}}}^{\theta_2}\,
\expt{-\frac{\alpha\,4^j\,r^2}{|z|}}\,
\frac{r^{2\,m}}{|z|^{m+3/2}}\,  {|dz|} \,
\Big(\aver{C_j(B)}
|h_{j}|^{p_0}\,dx\Big)^{\frac1{p_0}}
\\
&\lesssim&
2^{j\,(\theta_1- 2\,m-1)} \, \sum_{l=1}^{j} 2^{l}\,
\Big(\aver{2^{l+1}\,B}
|\nabla f|^{p_0}\,dx\Big)^{\frac1{p_0}},
\end{eqnarray*}
provided $2\, m+1>\theta_2$, where the last inequality follows by
repeating the calculations made to derive \eqref{eq:poincare}. The
term  corresponding to $\Gamma_-$ is controlled similarly. Plugging
both estimates into the representation of $\varphi(L)$ given by
\eqref{phi-L} one obtains
$$
\Big(\aver{B} |\varphi(L)
h_j|^{p_0}\,dx\Big)^{\frac{1}{p_0}}
\lesssim
2^{j\,(\theta_1- 2\,m-1)} \, \sum_{l=1}^{j} 2^{l}\,
\Big(\aver{2^{l+1}\,B}
|\nabla f|^{p_0}\,dx\Big)^{\frac1{p_0}}.
$$
The treatment for the term $j=2$ is similar using
$$|h_2|\le |f-f_{4\,B}|\,\bigchi_{8\,B\setminus 2\,B} \le  |f-f_{2\,B}|\,\bigchi_{8\,B\setminus 2\,B}
+  |f_{4\, B}-f_{2\,B}|\,\bigchi_{8\,B\setminus 2\,B}.
$$
Applying Minkowski's inequality and \eqref{eq:rep}, we obtain
\eqref{Tep-1}.  The lemma is proved.
\end{proof}

\noindent \textit{Case $p\in(\widetilde{p}_-,\widehat{p}_+)$}: Take
$p_0, q_0$ such that $\widetilde{p}_-<p_0<\widetilde{p}_+$ and
$p_0<p<q_0<\widehat{p}_+$.  Observe that $p_0\in \W_w\big(p_-,
p_+\big) \subset \Int\J_w(L)$ and $q_0\in \J_w(L)$. The proof
of Lemma \ref{lemma:est-reverseRiesz} extends \textit{mutatis
mutandis}  with $dw$ replacing $dx$ since there is an
$L^{p_0}$-Poincar\'e inequality for $dw$ (see Section \ref{sec:RT}). It
suffices to apply Theorem \ref{theor:main-w}  with underlying measure
$dw$ and no weight.   We leave further details to the reader.

\

\noindent \textit{Case $p\in(\max\big\{r_{w}\, ,\, (\widehat p_{-})_{w,*}\big\},\widetilde{p}_+)$}: It follows  a method in the
unweighted case by \cite{Aus} using an adapted Calder\'on-Zygmund
decomposition.

\begin{lemma}\label{lemmaCZD-w} Let $n\ge 1$, $w\in A_{\infty}$ and $1\le p< \infty$ such that $w\in A_{p}$. Assume that  $f\in \cals$ is  such that
$\|\nabla f\|_{L^p(w)} <\infty.$
Let $\alpha>0$. Then, one
can find a collection of balls $\{B_i\}_i$, smooth functions
$\{b_i\}_i$ and a function $g\in L^1_{\rm loc}(w)$ such that
\begin{equation}\label{eqcsds1}
f= g+\sum_i b_i
\end{equation} and the following properties hold:
\begin{equation}\label{eqcsds2}
|\nabla g(x)| \le C\alpha,
\quad
\text{for $\mu$-a.e. } x
\end{equation}
\begin{equation}\label{eqcsds3}
\supp b_i \subset B_{i}
\quad\text{and}\quad
\int_{B_i} |\nabla
b_i|^p\, dw \le C\alpha^p w(B_i),
\end{equation}
 \begin{equation}\label{eqcsds4}
\sum_i w(B_i) \le \frac{C}{\alpha^{p}} \int_{\re^n} |\nabla f|^p\, dw ,
\end{equation}
\begin{equation}\label{eqcsds5}
\sum_i \bigchi_{B_i} \le N,
\end{equation} where $C$ and  $N$ depends
only on the dimension, the doubling constant of $\mu$ and $p$. In
addition, for  $1\le q< p_{w}^*$, we have
 \begin{equation}\label{eqcsds6}
\Big(\aver{B_i} |b_i|^q\,dw\Big)^\frac1q
\lesssim
\alpha\, r(B_{i}).
\end{equation}
 \end{lemma}

\begin{proof}
Since $w\in A_p$, we have an $L^p(w)$ Poincar\'e inequality (see
\cite{FPW}). On the other hand, as $w\in A_p$ and $1\le q<p_w^*$ (if
$p=1$, i.e. $w\in A_{1}$, it holds also at $q=1_{w}^*= \frac
{n}{n-1}$ when $n\ge 2$) we can apply \cite[Corollary 3.2]{FPW} (when
checking the ``balance condition'' in that reference we have used
that $w\in A_r$ implies  $(|E|/|B|)^r\lesssim w(E)/w(B)$ for any ball
$B$ and any $E\subset B$). Thus there is an $L^p(w)-L^q(w)$ Poincar\'e
inequality:
\begin{equation}\label{Poin-p-q}
\Big(\aver{B}|f-f_{B,w}|^q \, dw\Big)^{\frac1q}
\lesssim \, r(B)\, \Big(\aver{B}|\nabla f|^p \, dw\Big)^{\frac1p}
\end{equation}
for all  locally Lipschitz functions $f$ and all balls $B$. These are
all the ingredients needed to invoke \cite[Proposition 9.1]{AM1}.
\end{proof}

We use the following resolution of $L^{1/2}$:
$$
L^{1/2}f= \frac 1 {\sqrt \pi}\,\int_0^\infty L e^{-t\,  L} f \,
\frac{dt}{\sqrt t}.
$$
  It suffices to work with $\int_\ep^R\ldots$,
to obtain bounds independent of $\ep,R$, and then to let
$\ep\downarrow 0$ and $R\uparrow \infty$: indeed, the truncated
integrals converge to $L^{1/2}f$ in $L^2$ when $f\in \cals$
and a use of Fatou's lemma concludes the proof.   For the truncated
integrals, all the calculations are justified. We write  $L^{1/2}$
where it is understood that it should be replaced by its
approximation at all places.

Take $q_0$ so that $\widetilde{p}_-<q_0<\widetilde{p}_+$. By the first case
of the proof,
\begin{equation}\label{eq:strongLq}
\|L^{1/2}f\|_{L^{q_0}(w)} \lesssim \|\nabla f\|_{L^{q_0}(w)}.
\end{equation}
We may assume that $\max\{\,r_{w}\, ,\, (\widehat
p_{-})_{w,*}\,\}<p<\widetilde p_{-}$, otherwise there is nothing to
prove. We claim that it is enough to show that
\begin{equation}
\label{eq:wtp-pReverse} \|L^{1/2}f\|_{L^{p,\infty}(w)} \lesssim
\|\nabla f\|_{L^p(w)}.
\end{equation}
Assuming this estimate we want to interpolate. To this end, we use the following lemma.

\begin{lemma}\label{lemma:density}
Assume $r>r_{w}$. Then $\D=\big\{(-\Delta)^{1/2}f\, : \, f\in \cals,
\supp \widehat f \subset \RR^n\setminus \{0\}\big\}$ is dense in
$L^r(w)$, where $\widehat f$ denotes  the Fourier transform of $f$.
\end{lemma}

\begin{proof} It is easy to see that $\D\subset \cals$ hence $\D\subset L^r(w)$.  As in \cite[p. 353]{Gra}, using that the classical Littlewood-Paley
series  converges in $L^r(w)$ since $w\in A_r$, it follows
that the set
$$
\widetilde{\D}=\big\{g\in \cals\, : \, \supp \widehat g \subset
\RR^n\setminus
\{0\},\ \supp \widehat g\, \mbox{ is compact}\big\}
$$
is dense in $L^r(w)$. We see that $\widetilde{\D}\subset \D$ and so $\D$
is dense in $L^r(w)$. For $g\in \widetilde{\D}$, $f=(-\Delta)^{-1/2}g$ is well-defined in $\cals$ as $\widehat{f}(\xi)=c\,|\xi|^{-1/2}\,\widehat{g}(\xi)$  and $\supp \widehat{f}\subset
\re^n\setminus\{0\}$. Hence, $g=(-\Delta)^{1/2}f \in \D$.
\end{proof}

If $r>r_{w}$,  the usual Riesz transforms, $\nabla (-\Delta)^{-1/2}$ , are bounded on $L^r(w)$  (this can be reobtained from the results in  Section  \ref{sec:RT}).   Also,  for $g\in L^r(w)$, one has
  $$
\|g\|_{L^r(w)} \sim \|\nabla (-\Delta)^{-1/2} g\|_{L^r(w)}
$$
using the identity $-I=R_1^2+\dots +R_n^2$ where
$R_j=\partial_j(-\Delta)^{-1/2}$. Thus,  for $g\in \D$,
$L^{1/2}(-\Delta)^{-1/2}g=L^{1/2}f$ if $f=(-\Delta)^{-1/2}g$ and $\|
\nabla  f\|_{L^r(w)} \sim   \|g \|_{L^r(w)}$  for  $r>r_{w}$. As
$r_w<p<q_0$,  \eqref{eq:strongLq} and \eqref{eq:wtp-pReverse}
reformulate into weighted strong type $(q_{0},q_{0})$ and weak type
$(p,p)$ of $T= L^{1/2}(-\Delta)^{-1/2}$ \textit{a priori} defined on
$\D$. Since $\D$ is dense in all $L^r(w)$  when $r>r_{w}$ by the
above lemma, we can extend $T$ by density in both cases and their
restrictions to the space of simple functions agree. Hence, we can
apply Marcinkiewicz interpolation and conclude again by density that
\eqref{eq:strongLq} holds for all $q$ with $p<q<q_0$ which leads to
the desired estimate.

\

Our goal is thus to establish  \eqref{eq:wtp-pReverse}, more precisely: for $f \in \cals$ and $\alpha>0$,
\begin{equation}\label{eq21}
w\{ |L^{1/2}f| >\alpha\} = w\{x \in \re^n: |L^{1/2}f(x)| >\alpha\} \le
\frac{C}{\alpha^p}\int_{\re^n} |\nabla f|^p\, dw.
\end{equation}
Since $p>r_{w}$, we have $w\in A_{p}$. From the condition $ (\widehat
p_{-})_{w,*}<p$, we have $\widehat p_{-} <p_{w}^*$. Therefore, there
exists $q\in (\widehat p_{-}, \widehat p_{+})=\Int\J_{w}(L)$ such
that $\widehat p_{-}<q<p_{w}^*$. Thus, we can apply the
Cal\-de\-r\'on-Zygmund decomposition of Lemma \ref{lemmaCZD-w} to $f$
at height $\alpha$ for  the measure $dw$ and write $f=g+\sum_i b_i$.
Using \eqref{eq:strongLq}, \eqref{eqcsds2} and $q_0>p$,  we have
\begin{align*}
w\Big\{|L^{1/2}g| >\frac \alpha 3 \Big\}
&\lesssim
\frac{1}{\alpha^{q_0}}\int_{\re^n} |L^{1/2} g|^{q_0}\, dw
\lesssim
\frac{1}{\alpha^{q_0}}\int_{\re^n} |\nabla g|^{q_0}\, dw
\lesssim
\frac{1}{\alpha^p}\int_{\re^n} |\nabla g|^p\,  dw
\\
&
\lesssim
\frac{1}{\alpha^p}\int_{\re^n} |\nabla f|^p\,  dw +
\frac{1}{\alpha^p}\int_{\re^n}
\Big|\sum_i\nabla b_i\Big|^p\,  dw
\lesssim
\frac{1}{\alpha^p}\int_{\re^n} |\nabla f|^p\,  dw,
\end{align*}
where the last estimate follows by applying \eqref{eqcsds5},
\eqref{eqcsds3}, \eqref{eqcsds4}.

To compute $L^{1/2}(\sum_i b_i)$, let $r_i=2^k$ if $2^k \le r(B_{i})
< 2^{k+1}$, hence $r_{i}\sim r(B_{i})$ for all $i$.  Write
$$
L^{1/2}
=
\frac 1 {\sqrt \pi}\,\int_0^{r_i^2} L e^{-t \, L}  \, \frac{dt}{\sqrt t}+ \frac 1 {\sqrt \pi}\,\int_{r_i^2}^\infty L e^{-t\,
L} \, \frac{dt}{\sqrt t}
=
T_i+U_i,
$$
and then
\begin{align*}
w\Big\{\Big|\sum_i L^{1/2} b_i\Big|>\frac{2\,\alpha}{3}\Big\}
&
\le
w\Big(\bigcup_i\, 4\,B_i\Big)+ w\Big\{\Big|\sum_i
U_ib_i\Big|>\frac{\alpha}{3}\Big\}
\\
&\hskip1.5cm
+ w\Big( \Big(\re^n\setminus\bigcup_i \, 4\,B_i \Big) \bigcap \Big\{
\Big|\sum_i T_ib_i\Big|
>\frac{\alpha}{3}\Big\}\Big)
\\
&
\lesssim
\frac1{\alpha^p}\,\int_{\re^n} |\nabla f|^p\,dw + I+II,
\end{align*}
where we have used \eqref{eqcsds4}. We estimate $II$. Since $q\in
\J_{w}(L)$, it follows that $t\,L\,e^{-t\, L}\in \offw{q}{q}$
by Proposition \ref{prop:sg-w:extension}, hence
\begin{align*}
II
&
\lesssim
\frac{1}{\alpha}\,\sum_i\sum_{j\ge 2} \int_{C_j(B_i)} |T_ib_i| \, dw
\lesssim
\frac{1}{\alpha}\,\sum_i\sum_{j\ge 2}
w(2^{j}\,B_i)\,\int_0^{r_i^2}\aver{C_{j}(B_{i})} |t\, L\, e^{-t \,L}
b_i |\, dw\,\frac{dt}{t^{3/2}}
\\
&\lesssim
\frac{1}{\alpha}\,\sum_i\sum_{j\ge 2} 2^{j\,D}\,w(B_i)\,\int_0^{r_i^2}
2^{j\,\theta_{1}}\, \dec{\frac{2^j\, r_i}{\sqrt t}}^{\theta_{2}} \,
\expt{-\frac{c\, 4^j\, r_{i}^2}{t}}\,\frac{dt}{t^{3/2}}
\, \Big(\aver{B_{i}} |b_{i}|^q\, dw\Big)^{\frac1q}
\\
&\lesssim
\sum_i\sum_{j\ge 2} 2^{j\,D}\,e^{-c\,4^j}\,w(B_i)
\lesssim
\sum_i w(B_i)
\lesssim \frac{1}{\alpha^p}\,\int_{\re^n}
|\nabla f|^p \, dw,
\end{align*}
where we have used  \eqref{eqcsds6} and \eqref{eqcsds4}, and $D$ is the doubling order of $dw$.

It remains to handling  the term $I$. Using functional calculus for
$L$ one can compute $U_i$ as $r_i^{-1}\psi(r_i^2L)$ with $\psi$ the
holomorphic function on the sector $\Sigma_{\pi /2}$
given by
\begin{equation}\label{eqpsiw}
\psi(z)= c \int_1^\infty z\, e^{-tz} \, \frac{dt}{\sqrt t}.
\end{equation}
It is easy to show that $|\psi(z)| \le C|z|^{1/2} e^{-c|z|}$,
uniformly on subsectors $\Sigma_{\mu}$, $0\le  \mu < {\frac \pi 2}$. We
claim that, since $q\in \Int\J_{w}(L)$,
 \begin{equation}\label{eq23w}
\Big\| \sum_{k\in \mathbb{Z}} \psi(4^kL)\ \beta_k \Big\|_{L^q(w)}
\lesssim
\bigg\|\Big(\sum_{k\in \mathbb{Z}}  |\beta_k|^2\Big)^{\frac12}
\bigg\|_{L^q(w)}.
\end{equation}
The proof of  this inequality is postponed until the end of Section
\ref{sec:sf}. We set $\beta_{k}= \sum_{i\, : \,
r_i=2^k}\frac{b_i}{r_i}$. Then,
$$
\sum_i U_i \, b_i
=
\sum_{k\in \mathbb{Z}} \psi(4^k\, L)
\bigg(\sum_{i\, : \, r_i=2^k}\frac{b_i}{r_i}\bigg)
=
\sum_{k\in \mathbb{Z}} \psi(4^k\,L) \beta_k.
$$
Using \eqref{eq23w}, the bounded overlap property \eqref{eqcsds5},
 \eqref{eqcsds6},   $r_{i}\sim r(B_{i})$ and \eqref{eqcsds4},
one has
\begin{align*}
I
&\lesssim
\frac1{\alpha^q}\,
\Big\|\sum_i U_i b_i\Big\|_{L^q(w)}^q
\lesssim
\frac1{\alpha^q}\,
\bigg\|\Big(\sum_{k\in \mathbb{Z}}  |\beta_k|^2\Big)^{\frac12}
\bigg\|_{L^q(w)}^q
\lesssim
\frac1{\alpha^q}\, \int_{\re^n} \sum_{i}  \frac{|b_i|^q}{r_i^q}\ dw
\\
&
\lesssim
\sum_{i}w(B_i)
\lesssim
\frac{1}{\alpha^p}\int_{\RR^n} |\nabla f|^p\,  dw.
\end{align*}
Collecting the obtained estimates, we conclude
\eqref{eq:wtp-pReverse} as desired.
\end{proof}

\begin{remark}\rm
If $w\in A_1$, $\W_w\big(p_-(L),p_+(L)\big)\ne \emptyset$ and
$(\widehat p_{-})_{w,*}<1$ then for all $f\in\cals$
$$
\|L^{1/2} f\|_{L^{1,\infty}(w)}
\lesssim
\|\nabla f\|_{L^1(w)}.
$$
This (that is \eqref{eq21} with $p=1$) uses a  similar argument (left
to the reader) once we have chosen an appropriate $q$ for which
$L^1(w)-L^q(w)$ Poincar\'e inequality holds:    since  $w\in A_1$, one
needs $q \le \frac{n}{n-1}$.  As $r_w=1$,  the assumption $(\widehat
p_{-})_{w,*}<1$ means that $\widehat{p}_{-} < \frac n {n-1}$ and so
we pick $q\in \Int \J_w(L)$ with $\widehat{p}_-<q<\frac{n}{n-1}$.
\end{remark}


\section{Square functions}\label{sec:sf}

We define the square functions for $x\in \RR^n$ and $f\in L^2$,
\begin{eqnarray*}
g_L f(x)
&=&
\Big(
\int_0^\infty |(t\,L)^{1/2}\,e^{-t\,L}f(x)|^2\,\frac{dt}{t}
\Big)^{\frac12},
\\[0.1cm]
G_L f(x)
&=&
\Big(
\int_0^\infty |\nabla e^{-t\,L}f(x)|^2\,dt
\Big)^{\frac12}.
\end{eqnarray*}
They are representative of a larger class of square functions and we
restrict our discussion to them to show the applicability of our methods.   They satisfy the following $L^p$
estimates.

\begin{theor}[\cite{Aus}]\label{theor:Aus-square}
$$
\Int \big\{ 1<p<\infty: \|g_L f\|_{p}\sim \|f\|_{p}, \forall\,f\in
L^p\cap L^2\big\}= \big(\,p_-(L),p_+(L)\big)
$$
and
$$
\Int \big\{1<p<\infty: \|G_L f\|_{p}\sim \|f\|_{p}, \forall\,f\in
L^p\cap L^2\big\}= \big(\,q_-(L),q_+(L)\big).
$$
\end{theor}

In this statement, $\sim$ can be replaced by $\lesssim$: the square
function estimates for $L$ (with $\lesssim$) automatically imply the reverse
ones for $L^*$. The part concerning $g_{L}$ can be obtained using an abstract
result of  Le Merdy \cite{LeM} as a consequence of the bounded
holomorphic functional calculus on $L^p$.  The method in \cite{Aus}
is direct. We remind the reader that in \cite{Ste}, these inequalities for $L=-\Delta$ were proved differently and the boundedness  of $G_{-\Delta}$ follows from that of $g_{-\Delta}$ and of the Riesz transforms $\partial_{j}(-\Delta)^{-1/2}$ (or vice-versa) using the commutation between $\partial_{j}$ and $e^{-t\, \Delta}$. Here, no such thing is possible.

We have the following weighted  estimates for square functions.

\begin{theor}\label{theor:square:weights}
Let $w\in A_\infty$.
\begin{list}{$(\theenumi)$}{\usecounter{enumi}\leftmargin=.8cm
\labelwidth=0.7cm\itemsep=0.3cm\topsep=.3cm
\renewcommand{\theenumi}{\alph{enumi}}}

\item If $\W_w\big(\,p_-(L),p_+(L)\big)\neq \emptyset$ and $p\in \Int
\J_w(L)$ then  for
all $f\in L^\infty_{c}$ we have
$$
\|g_L f\|_{L^p(w)}\lesssim \|f\|_{L^p(w)}.
$$

\item If $\W_w\big(\,q_-(L),q_+(L)\big)\neq \emptyset$ and $p\in \Int
\K_w(L)$ then  for
all $f\in L^\infty_{c}$ we have
$$
\|G_L f\|_{L^p(w)}\lesssim \|f\|_{L^p(w)}.
$$
\end{list}
\end{theor}

Note that the operators $(t\,L)^{1/2}\,e^{-t\,L}$ and $\nabla e^{-t\,L}$ extend to $L^p(w)$ when
 $p\in \Int\J_{w}(L)$ and $p\in \Int\K_{w}(L)$ respectively. By seeing $g_{L}$ and $G_{L}$ as  linear operators from scalar functions to $\HH$-valued functions (see below for  definitions), the above inequalities extend to all $f\in L^p(w)$ by density  (see the proof).

We also get reverse weighted square function estimates as follows.

\begin{theor}\label{theor:reverse-square:weights}
Let $w\in A_\infty$.

\begin{list}{$(\theenumi)$}{\usecounter{enumi}\leftmargin=.8cm
\labelwidth=0.7cm\itemsep=0.3cm\topsep=.3cm
\renewcommand{\theenumi}{\alph{enumi}}}

\item If $ \W_w\big(\,p_-(L),p_+(L)\big)\neq \emptyset$ and $p\in \Int\J_{w}(L)$ then
$$
\|f\|_{L^p(w)}\lesssim\|g_L f\|_{L^p(w)},
\qquad f\in L^p(w)\cap L^2.
$$

\item If  $r_{w}<p<\infty$,
$$
\|f\|_{L^p(w)}
\lesssim
\|G_L f\|_{L^p(w)},
\qquad f\in L^p(w)\cap L^2.
$$
\end{list}

\end{theor}

The restriction that $f\in L^2$
can be removed provided $g_{L}$ and $G_{L}$ are appropriately
interpreted: see the proofs.  We add a comment about sharpness of the
ranges of $p$ at the end of the section.

 As a corollary,  $g_{L}$ (resp. $G_{L}$) defines a new norm on $L^p(w)$ when  $p\in \Int\J_{w}(L)$ (resp. $p\in \Int\K_{w}(L)$ and $p > r_{w}$). Again,   Le Merdy's result cited above  \cite{LeM}  also gives such a result for $g_{L}$, but not for $G_{L}$. The restriction $p>r_{w}$ in part $(b)$ comes from the argument. We
do not know whether it is necessary for a given weight non identically 1.

Before we begin the arguments, we recall  some basic facts about
Hilbert-valued extensions of scalar inequalities. To do so we
introduce some notation: by $\HH$ we mean $L^2((0,\infty),\frac
{dt}{t})$ and $\normH{\cdot}$ denotes the norm in $\HH$. Hence, for a
function $h\colon \RR^n\times (0,\infty)\to \co$, we have for $x\in \RR^n$
$$
\normH{h(x,\cdot)}
=
\Big(\int_{0}^\infty |h(x,t)|^2\,\frac{dt}{t}\Big)^{1/2}.
$$
In particular, $$
g_{L}f(x)=\normH{\varphi(L,\cdot )f(x)}$$
with $\varphi(z,t)= (t\,z)^{1/2}\,e^{-t\,z}$
and
$$G_{L}f(x)= \normH{\nabla \varphi(L,\cdot )f(x)}$$
with $\varphi(z,t)= \sqrt t\,e^{-t\,z}$.
Let  $L^p_{\HH}(w)$ be the space of $\HH$-valued $L^p(w)$-functions equipped with the norm
$$
\|h\|_{L^p_{\HH}(w)}
=
\left(\int_{\re^n} \normH{h(x,\cdot)}^p\,dw(x)\right)^\frac1p.
$$

\begin{lemma}\label{lemma:M-Z:two-ops:new}
Let $\mu$ be a Borel measure on $\RR^n$ \textup{(}for instance, given
by an $A_{\infty}$ weight\textup{)}. Let    $1\le p \le q< \infty$.
Let $\D$ be a subspace of $\M$, the space of measurable functions in
$\RR^n$. Let
 $S,T $ be  linear operators from  $\D$ into $\M$.  Assume there exists $C_{0}>0$ such that  for all $f\in \D$, we have
$$
\|T f\|_{L^q(\mu)}
\le
C_0\, \sum_{j\ge 1} \alpha_{j }\|S f\|_{L^p(F_{j},\mu)},
$$
where $F_{j}$ are subsets of $\RR^n$ and $\alpha_{j}\ge 0$. Then,
there is an $\HH$-valued extension  with the same constant: for all
$f\colon \RR^n\times (0,\infty)\to \co$ such that  for
\textup{(}almost\textup{)} all $t>0, f(\cdot, t) \in \D$,
$$
\|T f\|_{L^q_{\HH}(\mu)}
\le
C_0\, \sum_{j\ge 1} \alpha_{j }\|S f\|_{L^p_{\HH}(F_{j},\mu)}.
$$
\end{lemma}

The extension of a linear operator  $T$ on $\co$-valued functions to
$\HH$-valued  functions is defined  for $x\in \RR^n$ and $t>0$ by $(T
h)(x,t)= T\big( h(\cdot,t)\big)(x)$, that is, $t$ can be considered
as a parameter and $T$ acts only on the variable in $\RR^n$. This
result is essentially the same as  the Marcinkiewicz-Zygmund theorem
and the fact that $\HH$ is isometric to $\ell^2$. That the norm
decreases uses $p\le q$. We refer to, for instance, \cite[Theorem
4.5.1]{Gra} for an argument that extends straightforwardly to our
setting.

\begin{proof}[Proof of Theorem \ref{theor:square:weights}. Part $(a)$] We split the argument in three  cases:
$p\in(\widetilde{p}_-,\widetilde{p}_+)$,
$p\in(\widetilde{p}_-,\widehat{p}_+)$,
$p\in(\widehat{p}_-,\widetilde{p}_+)$.

\

\noindent \textit{Case $p\in(\widetilde{p}_-,\widetilde{p}_+)$}:  By Proposition \ref{prop:weights}, there exist
$p_0, q_0$ such that
$$
p_-<p_0<p<q_0<p_+
\qquad {\rm and} \qquad
w\in A_{\frac{p}{p_0}}\cap RH_{\left(\frac{q_0}{p}\right)'}.
$$
We are going to apply Theorem \ref{theor:main-w} with $T=g_L$, $S=I$,
$\A_r=I-(I-e^{-r^2\,L})^m$,  $m$ large enough, underlying measure
$dx$ and weight $w$. We first see that \eqref{T:A} holds for all
$f\in L^\infty_{c}$. Here, we could have used the approach in
\cite{Aus}, but the one below adapts to the other two cases with
minor changes.

As $p_0, q_0\in \J(L)$ and $p_{0}\le q_{0}$, we know that $e^{-t\,L}\in\off{p_0}{q_0}$. If $B$ is a ball,   $j\ge 1$ and $g\in L^{p_{0}}$ with  $\supp g\subset C_j(B)$ we have
\begin{equation}\label{off-exten}
\Big( \aver{B} |e^{-k\,r^2\,L}g|^{q_0}\,dx
\Big)^{\frac1{q_0}}
\le
C_0\, 2^{j\,(\theta_1+\theta_2)}\,e^{-\alpha\,4^j}\,
\Big(
\aver{C_j(B)} |g|^{p_0}\,dx
\Big)^{\frac{1}{p_0}}.
\end{equation}
Lemma \ref{lemma:M-Z:two-ops:new} applied to $S=I$, $T\colon
L^{p_{0}}= L^{p_{0}}(\re^n, dx) \longrightarrow
L^{q_{0}}=L^{q_{0}}(\re^n, dx)$   given by
$$
T g= \Big(C_0\,
2^{j\,(\theta_1+\theta_2)}\,e^{-\alpha\,4^j}\Big)^{-1}
\,\frac{|2^{j+1}\,B|^{\frac1{p_0}}}{|B|^{\frac1{q_0}}}\, \bigchi_B \, e^{-k\,r^2\,L}
(\bigchi_{C_j(B)}\, g)
$$
yields
\begin{equation}\label{off-v-v}
\Big( \aver{B} \normH{e^{-k\,r^2\,L}g(x,\cdot)}^{q_0}\,dx
\Big)^{\frac1{q_0}}
\le C_0\, 2^{j\,(\theta_1+\theta_2)}\,e^{-\alpha\,4^j}\,
\Big(
\aver{C_j(B)} \normH{g(x,\cdot)}^{p_0}\,dx
\Big)^{\frac{1}{p_0}}
\end{equation}
for all $g\in L_{\HH}^{p_{0}}$ with  $\supp g(\cdot,t)\subset C_j(B)$ for each $t>0$.

As in \eqref{decomp-h}, for   $h \in L_{\HH}^{p_{0}}$ write $$h(x,t)=\sum_{j\ge
1}h_j(x,t), \quad x\in \re^n, \ t>0,
$$ where $h_j(x,t)=h(x,t)\,\bigchi_{C_j(B)}(x)$.  Using  \eqref{off-v-v}, we have for $1\le
k\le m$,
\begin{eqnarray}
\lefteqn{\hskip-2cm
\Big( \aver{B} \normH{e^{-k\,r^2\,L}h(x,\cdot)}^{q_0}\,dx
\Big)^{\frac1{q_0}}
\le
\sum_j
\Big(
\aver{B} \normH{e^{-k\,r^2\,L}h_j(x,\cdot)}^{q_0}\,dx
\Big)^{\frac1{q_0}}}\nonumber
\\
&\lesssim&
\sum_{j\ge 1} 2^{j\,(\theta_1+\theta_2)}\,e^{-\alpha\,4^j}\,
\Big(
\aver{2^{j+1}\,B} \normH{h(x,\cdot)}^{p_0}\,dx
\Big)^{\frac{1}{p_0}}.\label{estimate:q0-p0:ext}
\end{eqnarray}
Take $h(x,t)=(t\,L)^{1/2}\,e^{-t\,L} f(x)$.  Since $g_L
f(x)=\normH{h(x,\cdot)}$  and $f\in L^\infty_{c}$,  $h\in L_{\HH}^{p_{0}}$ by Theorem \ref{theor:Aus-square}  and
$$
g_L (e^{-k\,r^2\,L}f)(x)
=
\left(
\int_0^\infty |(t\,L)^{1/2}\,e^{-t\,L}\,e^{-k\,r^2\,L}f(x)|^2\,
\frac{dt}{t} \right)^{\frac12}
=
\normH{e^{-k\,r^2\, L}h(x,\cdot)}.
$$
Thus \eqref{estimate:q0-p0:ext} implies
$$
 \Big(\aver{B} |g_L
(e^{-k\,r^2\,L}f)|^{q_0}\,dx\Big)^{\frac1{q_0}}
\lesssim
\sum_{j\ge 1} 2^{j\,(\theta_1+\theta_2)}\,e^{-\alpha\,4^j}\,
\Big(
\aver{2^{j+1}\,B} |g_L f|^{p_0}\,dx
\Big)^{\frac{1}{p_0}}
$$
and it follows that $g_L$ satisfies \eqref{T:A}.

It remains to show that  \eqref{T:I-A}  with $Sf=f$ holds for all $f\in L^\infty_{c}$.
Write $f=\sum_{j\ge 1} f_j$ as
before. If $j=1$ we use that both $g_L$ and $(I-e^{-r^2\,L})^m$ are
bounded on $L^{p_0}$ (see Theorem \ref{theor:Aus-square} and
Proposition \ref{prop:sgfull}):
\begin{equation}\label{gL:f1}
\Big( \aver{B} |g_L(I-e^{-r^2\,L})^m f_1|^{p_0}\,dx\Big)^{\frac1{p_0}}
\lesssim
\Big( \aver{4\,B} |f|^{p_0}\,dx\Big)^{\frac1{p_0}}.
\end{equation}
For $j\ge 2$, we observe that
$$
g_L(I-e^{-r^2\,L})^m f_j(x)
=
\Big(
\int_0^\infty |(t\,L)^{1/2}e^{-t\,L}\,(I-e^{-r^2\,L})^mf_j(x)|^2
\,\frac{dt}{t} \Big)^{\frac12}
=
\normH{\varphi(L,\cdot)f_j(x)}
$$
where $\varphi(z,t)=(t\,z)^{1/2}\,e^{-t\,z}\,(1-e^{-r^2\,z})^m$.  As
in \cite{Aus}, the functions $\eta_{\pm}(\cdot,t)$ associated with  $\varphi(\cdot, t)$  by
\eqref{phi-L:eta}   verify
\begin{equation*}
|\eta_{\pm}(z,t)|
\lesssim \frac{t^{1/2}}{(|z|+t)^{3/2}}\,
\frac{r^{2\,m}}{(|z|+t)^m},
\qquad
z\in\Gamma_{\pm}, \ t >0.
\end{equation*}
Thus,
\begin{equation}\label{eq:gLeta}
\normH{\eta_{\pm}(z,\cdot)}
\le
\Big(
\int_{0}^\infty \frac{t}{(|z|+t)^{3}}\,
\frac{r^{4\,m}}{(|z|+t)^{2\,m}}\,\frac{dt}{t}\,\Big)^{\frac12}
\lesssim
\frac{r^{2\,m}}{|z|^{m+1}}.
\end{equation}
Next, applying Minkowski's inequality and $e^{-z\,L}\in
\off{p_0}{p_0}$, since $p_0\in\J(L)$, we have
\begin{eqnarray*}
\lefteqn{\hskip-1cm
\Big(\aver{B}
\bignormH{ \int_{\Gamma_+} e^{-z\,L} f_j(x)\,
\eta_+(z,\cdot)dz}^{p_0}\,dx\Big)^{\frac1{p_0}}}
\\
&\le&
\Big(\aver{B}
\Big( \int_{\Gamma_+} |e^{-z\,L}
f_j(x)|\,\normH{\eta_+(z,\cdot)}\,|dz|
\Big)^{p_0}\,dx\Big)^{\frac1{p_0}}
\\
&\le& \int_{\Gamma_+}
\Big(
\aver{B}|e^{-z\,L} f_j|^{p_0}\,dx
\Big)^{\frac1{p_0}}\,
\frac{r^{2\,m}}{|z|^{m+1}}\,|dz|
\\
&\lesssim&
2^{j\,\theta_1} \int_0^\infty
\dec{\frac{2^j\,r}{\sqrt{s}}}^{\theta_2}\,
\expt{-\frac{\alpha\,4^j\,r^2}{s}}\,
\frac{r^{2\,m}}{s^{m}}\,\frac{ds}{s} \,
\Big(
\aver{C_j(B)}|f|^{p_0}\,dx
\Big)^{\frac1{p_0}}
\\
&\lesssim&
2^{j\,(\theta_1-2\,m)}\, \Big(
\aver{C_j(B)}|f|^{p_0}\,dx
\Big)^{\frac1{p_0}}
\end{eqnarray*}
 provided $2\, m>\theta_2$. This plus the
corresponding term for $\Gamma_-$ yield
\begin{equation}\label{eq:ref2}
\Big(\aver{B} |g_L (I-e^{-r^2\,L})^m
f_j|^{p_0}\,dx\Big)^{\frac1{p_0}}
\lesssim 2^{j\,(\theta_1-2\,m)}\, \Big(
\aver{C_j(B)}|f|^{p_0}\,dx
\Big)^{\frac1{p_0}}.
\end{equation}
Collecting the latter estimate and \eqref{gL:f1}, we obtain that
\eqref{T:I-A} holds whenever $2\,m>\max\{\theta_1,\theta_{2}\}$.

\

\noindent \textit{Case $p\in(\widetilde{p}_-,\widehat{p}_+)$}: Take
$p_0, q_0$ such that $\widetilde{p}_-<p_0<\widetilde{p}_+$ and
$p_0<p<q_0<\widehat{p}_+$. Let $\A_r=I-(I-e^{-r^2\,L})^m$ for some
$m\ge 1$ to be chosen later. Remark that by the previous case, $g_L$
is bounded in $L^{p_{0}}(w)$ and so does $\A_{r}$ by Proposition
\ref{prop:sg-w:extension}. We apply Theorem
\ref{theor:main-w}  to $T=g_{L}$ and $S=I$ with underlying   measure $dw$ and no
weight: it is enough to see that $g_L$ satisfies \eqref{T:I-A} and
\eqref{T:A} on $L^\infty_{c}$. But this follows by adapting the preceding argument
replacing everywhere $dx$ by $dw$ and observing that $e^{-z\, L} \in \offw{p_{0}}{q_{0}}$. We skip details.

\

\noindent \textit{Case $p\in(\widehat{p}_-,\widetilde{p}_+)$}: Take
$p_0, q_0$ such that $\widetilde{p}_-<q_0<\widetilde{p}_+$ and
$\widehat{p}_-<p_0<p<q_0$. Set $\A_r=I-(I-e^{-r^2\,L})^m$ for some
integer $m\ge 1$ to be chosen later.  Since $q_0\in (\widetilde
p_-,\widetilde p_+)$,  by the first case, $g_L$ is bounded on
$L^{q_{0}}(w)$ and so does $\A_{r}$ by Proposition
\ref{prop:sg-w:extension}.  By Theorem \ref{theor:SHT:small} with
underlying Borel doubling measure $dw$, it is enough to show
\eqref{SHT:small:T:I-A} and \eqref{SHT:small:A}. Fix a ball $B$, $f
\in L^\infty_{c}$  supported on $B$.

Observe that \eqref{SHT:small:A} follows directly from
\eqref{est:Ar:w} since  $p_0, q_0\in \J_w(L)$ and $p_{0}\le q_{0}$.
We turn to \eqref{SHT:small:T:I-A}. Assume $j\ge 2 $.
 The argument is the same as the one for \eqref{eq:ref2} by reversing the roles of $C_{j}(B)$ and $B$, and using $dw$ and $ e^{-z\,L}\in\offw{p_0}{p_0}$ (since $p_{0}\in \J_{w}(L)$)  instead of $dx$ and $ e^{-z\,L}\in\off{p_0}{p_0}$. We obtain
\begin{equation*}
\label{eq:ref3bis}
\Big(\aver{C_{j}(B)} |g_{L} (I-e^{-r^2\,L})^m
f|^{p_0}\,dw\Big)^{\frac{1}{p_0}}
\lesssim
2^{j\,(\theta_1 -2\,m)} \, \Big(\aver{B}
|f|^{p_0}\,dw\Big)^{\frac1{p_0}}
\end{equation*}
provided $2\, m >\theta_{2}$ and it remains to impose further $2\, m > \theta_{1}+D$  to conclude, where $D$ is the doubling order of $w$.
\end{proof}

\begin{proof}[Proof of Theorem \ref{theor:square:weights}. Part $(b)$]
We split the argument in three  cases:
$p\in(\widetilde{q}_-,\widetilde{q}_+)$,
$p\in(\widetilde{q}_-,\widehat{q}_+)$,
$p\in(\widehat{q}_-,\widetilde{q}_+)$.

\

\noindent \textit{Case $p\in(\widetilde{q}_-,\widetilde{q}_+)$}:  By Proposition \ref{prop:weights}, there exist
$p_0, q_0$ such that
$$
q_-<p_0<p<q_0<q_+
\qquad {\rm and} \qquad
w\in A_{\frac{p}{p_0}}\cap RH_{\left(\frac{q_0}{p}\right)'}.
$$
We are going to apply Theorem \ref{theor:main-w}  with underlying
measure $dx$ and weight $w$ to $T=G_L$, $S=I$,
$\A_r=I-(I-e^{-r^2\,L})^m$, $m$ large enough. We begin with
\eqref{T:A}. Fix $1\le k\le m$ and $B$ a ball.  Combining
\eqref{eq:est-Riesz:2} and Lemma \ref{lemma:M-Z:two-ops:new} with
$T=\nabla e^{-k\, r^2\, L}$ and $S=\nabla$, we obtain
\begin{equation*}
\Big(\aver{B} \normH{\nabla
e^{-k\,r^2\,L}h(x,\cdot)}^{q_0}\,dx\Big)^{\frac1{q_0}}
\lesssim
 \sum_{j\ge 1} g_{2}(j)\,\Big(
\aver{2^{j+1}\,B}\normH{\nabla h(x,\cdot)}^{p_0}\,dx
\Big)^{\frac1{p_0}}\end{equation*}
with $g_2(j)=C_m\,2^j\,\sum_{l\ge j} 2^{l\,\theta}\,e^{-\alpha\,4^l}$
for some $\theta>0$ whenever $h\colon \RR^n \times (0,\infty)
\longrightarrow \co$ is such that $h$ and $\nabla h$ belong to
$L^{p_{0}}$ (our space $\D$). Setting
$h(x,t)=\sqrt{t}\,e^{-t\,L}f(x)$ for $f\in L^\infty_{c}$, we note
that  $h(\cdot, t) \in L^{p_{0}}$ and  $\nabla h(\cdot, t) \in
L^{p_{0}}$ for each $t>0$. Hence,  the above estimate applies. Since
$\normH{\nabla h(x,\cdot)}=G_L f(x)$ and $\normH{\nabla
e^{-k\,r^2\,L}h(x,\cdot)}=G_L(e^{-k\,r^2\,L}f)(x)$, we obtain
$$
\Big(\aver{B} |G_L(e^{-k\,r^2\,L}f)|^{q_0}\,dx\Big)^{\frac1{q_0}}
\lesssim \sum_{j\ge 1} g_{2}(j)\,\Big(
\aver{2^{j+1}\,B}|G_L f|^{p_0}\,dx
\Big)^{\frac1{p_0}},
$$
which is \eqref{T:A} after expanding $\A_r$.

It remains to checking \eqref{T:I-A} for $G_{L}$ and $S=I$ for $f\in L^\infty_{c}$.
Fix a ball $B$. As before, write
$f=\sum_{j\ge 1} f_j$ where $f_j=f\,\bigchi_{C_j(B)}$. Since $p_0\in
\Int \K(L)$, both $G_L$ and $(I-e^{-r^2\,L})^m$ are bounded on
$L^{p_0}$ by Theorem \ref{theor:Aus-square}  and Proposition
\ref{prop:sgfull}. Then for $j=1$ we have
\begin{equation}\label{GL:f1}
\Big(
\aver{B} |G_L (I-e^{-r^2\,L})^m f_1|^{p_0}\,dx \Big)^{\frac1{p_0}}
\lesssim
\Big(
\aver{4\,B} |f|^{p_0}\,dx\Big)^{\frac1{p_0}}.
\end{equation}
For $j\ge 2$, we observe that
$$
G_L(I-e^{-r^2\,L})^m f_j(x)
=
\Big(
\int_0^\infty |\sqrt{t}\nabla e^{-t\,L}\,(I-e^{-r^2\,L})^mf_j(x)|^2
\frac{dt}{t} \Big)^{\frac12}
=
\normH{\nabla \varphi(L,\cdot)f_j(x)}
$$
where $\varphi(z,t)= \sqrt{t} \,e^{-t\,z}\,(1-e^{-r^2\,z})^m$. As
in \cite{Aus}, the functions $\eta_{\pm}(\cdot,t)$ associated for each $t>0$ with $\varphi(\cdot, t)$  by
\eqref{phi-L:eta}   verify
$$
|\eta_{\pm}(z, t)|
\lesssim \frac{\sqrt{t}}{|z|+t}\,  \frac{r^{2\,m}}{(|z|+t)^m},
\qquad
z\in\Gamma_{\pm}, \ t>0,
$$
and so
\begin{equation}\label{eta-GL-L2}
\normH{\eta_{\pm}(z,\cdot)}
\le
\Big(\int_{0}^\infty \frac{t}{(|z|+t)^{2}}\,
\frac{r^{4\,m}}{(|z|+t)^{2\,m}}\,\frac{dt}{t}\,\Big)^{\frac12}
\lesssim
\frac{r^{2\,m}}{|z|^{m+1/2}}.
\end{equation}
Using Minkowski's inequality and
 $\sqrt{z}\,\nabla e^{-z\,L}\in \off{p_0}{p_0}$ since $p_0\in \K(L)$,
\begin{eqnarray*}
\lefteqn{\hskip-1cm
\Big(\aver{B}
\bignormH{\int_{\Gamma_+} \nabla e^{-z\,L} f_j(x)\,
\eta_+(z,\cdot)\,dz }^{p_0}\,dx\Big)^{\frac1{p_0}}}
\\
&\le&
\Big(\aver{B}
\Big( \int_{\Gamma_+} |\sqrt{z}\,\nabla e^{-z\,L}
f_j(x)|\,\normH{\eta_+(z,\cdot)}\,\frac{|dz|}{|z|^{1/2}}
\Big)^{p_0}\,dx\Big)^{\frac1{p_0}}
\\
&\le& \int_{\Gamma_+}
\Big(
\aver{B}|\sqrt{z}\,\nabla  e^{-z\,L} f_j(x)|^{p_0}\,dx
\Big)^{\frac1{p_0}}\,
\frac{r^{2\,m}}{|z|^{m+1/2}}\,\frac{|dz|}{|z|^{1/2}}
\\
&\lesssim&
2^{j\,\theta_1} \int_0^\infty
\dec{\frac{2^j\,r}{\sqrt{s}}}^{\theta_2}\,
\expt{-\frac{\alpha\,4^j\,r^2}{s}}\,
\frac{r^{2\,m}}{s^{m}}\,\frac{ds}{s} \,
\Big(
\aver{C_j(B)}|f|^{p_0}\,dx
\Big)^{\frac1{p_0}}
\\
&\lesssim&
2^{j\,(\theta_1-2\,m)}\, \Big(
\aver{C_j(B)}|f|^{p_0}\,dx
\Big)^{\frac1{p_0}}
\end{eqnarray*}
 provided $2\, m>\theta_2$. This, plus the
corresponding term for $\Gamma_-$, yields
\begin{equation}\label{eq:ref3}
\Big(\aver{B} |G_L (I-e^{-r^2\,L})^m
f_j|^{p_0}\,dx\Big)^{\frac1{p_0}}
\lesssim
2^{j\,(\theta_1-2\,m)}\, \Big(
\aver{C_j(B)}|f|^{p_0}\,dx
\Big)^{\frac1{p_0}}.
\end{equation}
Collecting the latter estimate and \eqref{GL:f1}, we obtain by
Minkowski's inequality
\begin{eqnarray*}
\Big(\aver{B} |G_L (I-e^{-r^2\,L})^m
f|^{p_0}\,dx\Big)^{\frac1{p_0}}
&\lesssim&
\sum_{j\ge 1}2^{j\,(\theta_1-2\,m)}\, \Big(
\aver{C_j(B)}|f|^{p_0}\,dx
\Big)^{\frac1{p_0}}.
\end{eqnarray*}
 Therefore,   \eqref{T:I-A} holds on taking $2\, m >\sup(\theta_{1},\theta_{2}).$

\

\noindent \textit{Case $p\in(\widetilde{q}_-,\widehat{q}_+)$}: Take
$p_0, q_0$ such that $\widetilde{q}_-<p_0<\widetilde{q}_+$ and
$p_0<p<q_0<\widehat{q}_+$. Let $\A_r=I-(I-e^{-r^2\,L})^m$ for some
$m\ge 1$ to be chosen later. As   $p_0\in (\widetilde q_-, \widetilde
q_+)$,  both $G_L$ and $\A_{r}$  are bounded on $L^{p_0}(w)$ (we have
just shown it for $G_L$ and Proposition \ref{prop:sg-w:extension}
yields it for $\A_{r}$ with a uniform norm in $r$). By Theorem
\ref{theor:main-w} with underlying doubling measure $dw$ and no
weight, it is enough to verify \eqref{T:I-A} and \eqref{T:A}  on
$\D=L^\infty_{c}$ for $T=G_L$, $S=I$. It suffices to  copy the
preceding argument replacing everywhere $dx$ by $dw$,  observing that
$p_0, q_0\in \K_w(L)$ implies weighted off-diagonal estimates and an
$L^{p_0}(w)$ Poincar\'e inequality, and applying Lemma
\ref{lemma:M-Z:two-ops:new}  to obtain an $\HH$-valued extension. We
leave the details to the reader.

\

\noindent \textit{Case $p\in(\widehat{q}_-,\widetilde{q}_+)$}:
Take $p_0, q_0$ such that $\widetilde{q}_-<q_0<\widetilde{q}_+$ and
$\widehat{q}_-<p_0<p<q_0$. Set $\A_r=I-(I-e^{-r^2\,L})^m$ for some
$m\ge 1$ to be chosen later. Since $q_0\in
(\widetilde{q}_-,\widetilde{q}_+) $, it follows that $G_L$ is bounded
on $L^{q_0}(w)$ and so is $\A_r$  by Proposition
\ref{prop:sg-w:extension}.   By Theorem \ref{theor:SHT:small} with
underlying measure $dw$, it is enough to show \eqref{SHT:small:T:I-A}
and \eqref{SHT:small:A}.

Observe that \eqref{SHT:small:A} is nothing but  \eqref{est:Ar:w}
since $p_0, q_0\in \K_w(L) \subset \J_{w}(L)$. The proof of
\eqref{SHT:small:T:I-A} is again  analogous to \eqref{eq:ref3} in the
weighted setting exchanging the roles of $C_{j}(B)$ and $B$. We skip
details.
\end{proof}

To prove Theorem \ref{theor:reverse-square:weights},  part $(a)$, we introduce the following operator.
Define  for $f\in L^2_{\HH}$ and  $x\in \RR^n$,
\begin{eqnarray*}
T_L f(x)
=\int_0^\infty (t\,L)^{1/2}\,e^{-t\,L}f(x,t )\,\frac{dt}{t}.
\end{eqnarray*}
Recall that $(t\,L)^{1/2}\,e^{-t\,L}f(x,t
)=(t\,L)^{1/2}\,e^{-t\,L}(f(\cdot,t ))(x)$. Hence, $T_{L}$ maps
$\HH$-valued functions to $\co$-valued functions.  We note that, for
$f \in L^2_{\HH}$ and $h\in L^2$, we have
$$
\int_{\re^n} T_{L}f\, \overline h\, dx = \int_{\re^n} \int_0^\infty f(x,t) \, \overline{ (t\,L^*)^{1/2}\,e^{-t\,L^*}h(x )}\,\frac{dt}{t}dx,
$$
where $L^*$ is the adjoint (on $L^2$) of $L$, hence,
$$
\Big| \int_{\re^n} T_{L}f\, \overline h\, dx\Big| \le \int_{\re^n} \normH{f(x, \cdot)} \, g_{L^*}(h)(x) \, dx.
$$
Let $p_{-}(L)<p<p_{+}(L)$.  Since $p_-(L^*)=\big(p_+(L)\big)'<p'<\big(p_-(L)\big)'=p_+(L^*)$,  $g_{L^*}$ is bounded on $L^{p'}$. This and a density argument imply that $T_{L}$ has a bounded extension from $L^{p}_{\HH}$ to $L^{p}$. The weighted version is as follows.

\begin{theor}\label{theor:TL}
Let $w\in A_\infty$.
 If $\W_w\big(\,p_-(L),p_+(L)\big)\neq \emptyset$ and $p\in \Int
\J_w(L)$ then  for
all $f\in L^\infty_{c}(\re^n \times (0,\infty))$ we have
$$
\|T_L f\|_{L^p(w)}\lesssim \|f\|_{L^p_{\HH}(w)}.
$$
Hence, $T_L$ has a bounded extension from $L^p_{\HH}(w)$ to $L^p(w)$.
\end{theor}

The duality argument above
 works for exponents in
$\W_w\big(\,p_-(L),p_+(L)\big)$, but we do not know how to extend it
to all of  $\Int \J_w(L)$. Hence, we proceed via a direct proof where
duality is  used only when $w=1$.

\begin{proof} We split the argument in three  cases:
$p\in(\widetilde{p}_-,\widetilde{p}_+)$,
$p\in(\widetilde{p}_-,\widehat{p}_+)$,
$p\in(\widehat{p}_-,\widetilde{p}_+)$.

\

\noindent \textit{Case $p\in(\widetilde{p}_-,\widetilde{p}_+)$}: By
Proposition \ref{prop:weights}, there exist $p_0, q_0$ such that
$$
p_-<p_0<p<q_0<p_+
\qquad {\rm and} \qquad
w\in A_{\frac{p}{p_0}}\cap RH_{\left(\frac{q_0}{p}\right)'}.
$$
We
are going to apply Theorem \ref{theor:main-w} (in fact, its
vector-valued extension) with underlying measure $dx$ and weight $w$
to the linear operator $T=T_L$ with $S=I$ and
$\A_r=I-(I-e^{-r^2\,L})^m$, $m$ large enough. Here,  $\A_{r}$ denotes
 both the scalar operator and its $\HH$-valued extension. We first see that $T_L$ satisfies \eqref{T:A} with $p_{0},q_{0}$ for  $f\in L_{c}^\infty(\RR^n\times (0,\infty))$. Let
$B$ be a ball. Note
that $T_{L}\A_{r}f= \A_{r}T_{L}f$ with our confusion of notation. Hence \eqref{T:A} is a simple
consequence of \eqref{off-exten} applied to $g=T_{L}f$.

Next, it remains to check  \eqref{T:I-A}. Let $f \in
L_{c}^\infty(\RR^n\times (0,\infty))$ and  let $B$ be a ball. As in
\eqref{decomp-h}, we write
$$
f(x,t)=\sum_{j\ge 1}f_j(x,t),
$$
where $f_j(x,t)=f(x,t)\,\bigchi_{C_j(B)}(x)$.  For $T_{L}(I-\A_{r})
f_{1}$, we use the boundedness of $T_{L}$ from $L^{p_{0}}_{\HH}$ to
$L^{p_{0}}$ noted above and  the $L^{p_{0}}_{\HH}$ boundedness of
$\A_{r}$ to obtain
$$
\Big(\aver{B} |T_{L}(I-\A_{r}) f_{1}|^{p_{0}}\, dx\Big)^{\frac 1 {p_{0}}} \lesssim
 \Big(\aver{4\, B} \normH{f(x, \cdot)}^{p_{0}}\, dx\Big)^{\frac 1 {p_{0}}}.
$$
For $j\ge 2$, the functions $\eta_{\pm}(z,t)$ associated with
$\varphi(z,t)=(t\,z)^{1/2}\, e^{-t\,z}\,  (1-e^{-r^2\,z})^m$ by
\eqref{phi-L:eta} satisfy \eqref{eq:gLeta}. Hence,
\begin{eqnarray*}
\lefteqn{\hskip-1cm
\Big(
\aver{B} \left| \int_{0}^\infty \int_{\Gamma_{+}} e^{-z\,
L}f_{j}(x,t) \eta_{+}(z,t)\, dz\frac {dt}t \right |^{p_{0}}\,
dx
\Big)^{\frac 1 {p_{0}}}}
\\
&\lesssim&
\Big(
\aver{B} \Big(  \int_{\Gamma_{+}} \normH{e^{-z\, L}f_{j}(x,
\cdot)}\,\normH{ \eta_{+}(z,\cdot)}\, |dz| \Big)^{p_{0}}\, dx
\Big)^{\frac 1 {p_{0}}}
\\
&\lesssim&
\int_{\Gamma_{+}} \Big(\aver{B}   \normH{e^{-z\, L}f_{j}(x,
\cdot)}^{p_{0}}\, dx\Big)^{\frac 1 {p_{0}}} \normH{
\eta_{+}(z,\cdot)}\, |dz|
\\
&\lesssim&
2^{j\, \theta_{1}}\int_{0}^\infty \dec{\frac{2^j \, r}{\sqrt
s}}^{\theta_{2}} \, \expt{-\frac{c\, 4^j\, r^2}{s}} \, \frac
{r^{2m}}{s^m}\, \frac {ds}s \,
\Big(\aver{C_{j}(B)}   \normH{f(x, \cdot)}^{p_{0}}\, dx\Big)^{\frac 1 {p_{0}}}
\\
&\lesssim&
2^{j\, (\theta_{1}-2\,m)}\,
\Big(\aver{C_{j}(B)}   \normH{f(x, \cdot)}^{p_{0}}\, dx\Big)^{\frac 1 {p_{0}}}
\end{eqnarray*}
 where we used the $\HH$-valued extension of
$e^{-z\, L} \in \off{p_{0}}{p_{0}}$ and assumed $2\,m > \theta_{2}$.
This, plus the corresponding term for $\Gamma_-$, yields
\begin{equation}
\label{eq:ref4}
\Big(\aver{B} |T_{L}(I-\A_{r}) f_{j}|^{p_{0}}\, dx\Big)^{\frac 1{p_{0}}}
\lesssim
2^{j\, (\theta_{1}-2\,m)}\,
\Big(\aver{C_{j}(B)}   \normH{f(x, \cdot)}^{p_{0}}\, dx\Big)^{\frac 1 {p_{0}}}
\end{equation}
and therefore \eqref{T:I-A} follows on taking also $2\,m >
\theta_{1}$.

 \

\noindent \textit{Case $p\in(\widetilde{p}_-,\widehat{p}_+)$}: Take
$p_{0}, q_{0}$ with $\widetilde p_{-}< p_{0}<\widetilde p_{+}$
and $p_{0}<p<q_{0}<\widehat p_{+}$. It suffices to apply  the
$\HH$-valued extension of Theorem \ref{theor:main-w}  with
underlying doubling measure $dw$ and no weight to the linear operator
$T=T_L$, $S=I$ and $\A_r=I-(I-e^{-r^2\,L})^m$,  $m$ large enough.
This is done exactly as in the previous case. At some step we have to
use that $T_L$ is bounded from $L^{p_0}_{\HH}(w)$ to $L^{p_0}(w)$
which follows by the previous case. We leave the details to the
reader.

\

\noindent \textit{Case $p\in(\widehat{p}_-,\widetilde{p}_+)$}: Take
$p_{0}, q_{0}$ with $\widetilde p_{-}<q_{0}<\widetilde p_{+}$
and $\widehat p_{-}<p_{0}<p<q_{0}$.  Since $q_0\in (\widetilde
p_-,\widetilde p_+)$,  by the first case, $T_L$ is bounded from
$L^{q_{0}}_{\HH}(w)$ to $L^{q_{0}}(w)$ and so does
$\A_r=I-(I-e^{-r^2\,L})^m$,  $m\ge 1$, on $L^{q_{0}}_{\HH}(w)$ by
Proposition \ref{prop:sg-w:extension} and Lemma \ref{lemma:M-Z:two-ops:new}.
 By Theorem \ref{theor:SHT:small} (in fact, its $\HH$-valued extension)   with
underlying Borel doubling measure $dw$, it is enough to show
\eqref{SHT:small:T:I-A} and \eqref{SHT:small:A} on $\D=L^\infty_{c}(\re^n \times (0,\infty))$
 for $T=T_{L}$ and $\A_{r}$ with large enough $m$.  As usual, the
latter is a mere consequence of $e^{-t\, L} \in \offw{q_{0}}{q_{0}}$
and its $\HH$-valued analog. The first condition is again a repetition of the argument for
\eqref{eq:ref4} in the weighted setting switching $C_{j}(B)$ and $B$. We skip details.
\end{proof}

\begin{proof}[Proof of Theorem \ref{theor:reverse-square:weights}]
We begin with part $(a)$.  Fix $p\in\Int \J_w(L)$ where $w\in
A_\infty$ so that $\W_w\big(p_-(L),p_+(L)\big) \ne \emptyset$. Let
$f\in L^2$ and define $F$ by $F(x,t)=(t\, L)^{1/2} \, e^{-t\,
L}f(x)$. Note that $F\in  L^2_{\HH}$ since
$\|F\|_{L^2_{\HH}}=\|g_{L}f\|_{2}$. By functional calculus on $L^2$,
we have
\begin{equation}
\label{eq:repTL}
f= 2\int_{0}^\infty (t\, L)^{1/2} \, e^{-t\, L}F(\cdot,t)\, \frac{dt}t = 2 \, T_{L} F
\end{equation}
with convergence in $L^2$. Note that for $p\in\Int\J_{w}(L)$,
$e^{-t\, L}$ has an infinitesimal generator on $L^p(w)$ as recalled
in Remark \ref{remark:inf-gen}. Let us call $L_{p,w}$ this generator.
In particular $e^{-t\, L}$ and $e^{-t\, L_{p,w}}$ agree on $L^p(w)
\cap L^2$. Our results assert that $L_{p,w}$ has a bounded
holomorphic functional calculus on $L^p(w)$, hence replacing $L$ by
$L_{p,w}$ and $f\in L^2$ by $f\in L^p(w)$, we see that $F\in
L^p_{\HH}(w)$ with $\|F\|_{L^p(w)}=\|g_{L_{p,w}}f\|_{L^p(w)}$  and
\eqref{eq:repTL} is valid with convergence in $L^p(w)$ (this is
standard fact from functional calculus and we skip details). Thus, by
Theorem \ref{theor:TL},
$$
\|f\|_{L^p(w)} = 2 \|T_{L_{p,w}}F\|_{L^p(w)} \lesssim \|F\|_{L^p_{\HH}(w)} = \|g_{L_{p,w}}f\|_{L^p(w)}.
$$
Noting  that $g_{L}f= g_{L_{p,w}}f$ when $f\in L^2\cap L^p(w)$ and
$T_{L}F=T_{L_{p,w}}F$ when $F\in L^2_{\HH}\cap L^p_{\HH}(w)$, part
$(a)$ is proved.

\

Let us show part $(b)$, that is the corresponding inequality for $G_L$. Fix $w\in
A_\infty$. We use the following
estimate from \cite{Aus}: for $f,h\in L^2$
$$
\left|\int_{\re^n} f\,\overline{h}\,dx\right|
\le
(1+\|A\|_\infty)\, \int_{\re^n} G_Lf\,G_{-\Delta}h\,dx,
$$
where $G_{-\Delta}$ is the square function associated with the
operator $-\Delta$. It is well known that $G_{-\Delta}$ is bounded on
$L^q(u)$ for all $1<q<\infty$ and all $u\in A_q$. Let us emphasize
that, indeed, the results that we have proved can be applied to the
operator $-\Delta$ and so $G_{-\Delta}$ is bounded on $L^q(u)$ for
$u\in A_\infty$ and all $
q\in\W_{u}\big(q_-(-\Delta),q_+(-\Delta)\big)= \W_{u}(1,\infty)$,
that is, for all $1<q<\infty$ and  $u\in A_q$.

Coming back to the argument, let   $p>r_{w}$, hence $w\in A_{p}$.  Let $f\in L^2\cap L^p(w)$. Then
$$
\int_{\re^n} \vert f\vert ^p\,  dw= \lim_{N,k, R\to \infty} \int_{\re^n} f\,  \overline{h}\,  dw_{N}
$$
with  $w_N=\min\{w,N\}$ and $h= f \vert f\vert ^{p-2}
\bigchi_{B(0,R)}\bigchi_{\{0<\vert f\vert \le k\}}$. Note that
$\|h\|_{L^{p'}(w_N)}\le \|f\|_{L^p(w)}^{p-1}$ and that $hw_{N}$ is a
bounded compactly supported function, hence in $L^2$.

Observe that $w_{N}\in A_p$ with $A_p$-constant smaller than the one
for $w$.   As observed, $G_{-\Delta}$ is bounded on
$L^{p'}(w_{N}^{1-p'})$ since $w_{N}^{1-p'}\in A_{p'}$. Thus,
 we have
\begin{align*}
\Big|\int_{\re^n} f\,\overline{h}\,dw_{N} \Big|
&=
\Big|\int_{\re^n} f\,\overline{hw_{N}}\,dx \Big|
\le
(1+\|A\|_\infty)\, \int_{\re^n} G_L f\,G_{-\Delta}(h\,w_{N})\,dx
\\
&\le
(1+\|A\|_\infty)\, \|G_L
f\|_{L^p(w_{N})}\,\|G_{-\Delta}(h\,w_{N})\|_{L^{p'}(w_{N}^{1-p'})}
\\
&\le
C\,\|G_L f\|_{L^p(w_{N})}\,\|h\,w_{N}\|_{L^{p'}(w_{N}^{1-p'})}
\\
&\le
C\, \|G_L f\|_{L^p(w)} \, \| f\|_{L^p(w)}^{p-1}
\end{align*}
with $C$ is independent of $N,k,R$  and where  we have used that $w_{N}\le w$. Thus taking limits
$N\to\infty$ first and then $k\to \infty$ and $R\to \infty$, we obtain
$$
\| f\|_{L^p(w)}^{p}\le C \|G_L f\|_{L^p(w)} \, \| f\|_{L^p(w)}^{p-1}.
$$
\end{proof}

\begin{proof}[Proof of  \eqref{eq23w}]  The operator in \eqref{eq23w} is similar to $T_{L}$,
 changing continuous times $t$ to discrete
times $4^k$ and   $z^{1/2}e^{-z}$ to $\psi(z)$.  Since $\psi(z)$ has the same  quantitative
properties as $z^{1/2}e^{-z}$   (decay at 0 and at
infinity), the  proof of Theorem \ref{theor:TL}  applies and
furnishes \eqref{eq23w}.
\end{proof}

\begin{remark}
\rm  $\Int\J_{w}(L)$ is the sharp range  up to endpoints for
$\|g_{L}f\|_{L^p(w)} \sim \|f\|_{L^p(w)}$. Indeed, we have $g_{L}(e^{- t\, L} f) \le g_{L}f$  for all
$t>0$. Hence,  the equivalence implies the uniform
$L^p(w)$ boundedness of $e^{-t\, L}$, which implies $p\in
\widetilde\J_{w}(L)$ (see Proposition \ref{prop:sg-w:extension}). Actually,  $\Int\J_{w}(L)$ is also the sharp range up to endpoints for the inequality
$\|g_{L}f\|_{L^p(w)} \lesssim \|f\|_{L^p(w)}$. It suffices to adapt the interpolation procedure in \cite[Theorem 7.1, Step 7]{Aus}. We skip details.

Similarly, this interpolation procedure also shows that  $\Int\K_{w}(L)$ is also  sharp up to endpoints for  $\|G_{L}f\|_{L^p(w)} \lesssim \|f\|_{L^p(w)}$.
\end{remark}

\section{Some vector-valued estimates}\label{sec:vv}

In \cite{AM1}, we also obtained vector-valued inequalities.

\begin{prop}
Let $\mu, p_{0},q_{0}, T,\A_{r}, \D$ be as in   Theorem
\ref{theor:main-w}  and assume \eqref{T:I-A} and \eqref{T:A} with
$S=I$.  Let $p_{0}<p,r<q_{0}$. Then, there is a constant $C$ such
that for all $f_k\in \D$
\begin{equation}\label{T:v-v}
\Big\|
\Big(\sum_k |T f_k|^r\Big)^{\frac1r}
\Big\|_{L^p(\mu)}
\le
C\,\Big\|
\Big(\sum_k | f_k|^r\Big)^{\frac1r}
\Big\|_{L^p(\mu)}.
\end{equation}
\end{prop}

Let  us see how it applies here.

First, let $T=\varphi(L)$ ($\varphi$ bounded holomorphic in an appropriate
sector).  Theorem \ref{theor:B-K:weights} says that $T$ is bounded on
$L^p(w)$ for all $p \in \Int\J_{w}(L)$. Also, for $p_{0},q_{0}\in \Int\J_{w}(L)$ with $p_{0}<q_{0}$,
we have $L^{p_{0}}(w)-L^{q_{0}}(w)$ off-diagonal
estimates on balls for $\A_{r}= I- (I-e^{-r^2L})^m$. Hence,
we can prove  \eqref{T:I-A} and \eqref{T:A} with
$S=I$ where $dx$  is now replaced by $w(x) dx$ by mimicking   the first
case of the proof  of Theorem \ref{theor:B-K:weights} in the weighted context.  Hence, one can apply the proposition above with $d\mu=w \, dx$ to above weighted vector-valued estimates for $\varphi(L)$  with all $p,r \in \Int\J_{w}(L)$.

The same  weighted vector-valued estimates hold  with all $p,r \in \Int\J_{w}(L)$  with $T=g_{L}$ starting  from  Theorem \ref{theor:square:weights} and mimicking the proof of its first case with $dx$ replaced with $w(x) dx$.

If $T=\nabla L^{-1/2}$ or $T=G_{L}$, then the same reasoning  applies modulo  the
Poincar\'e inequality used towards obtaining \eqref{T:A}. Hence, we conclude that for
both $\nabla L^{-1/2}$ and $G_{L}$, one has  \eqref{T:v-v}  with
$d\mu= w dx$ and $p,r\in \Int\K_{w}(L) \cap (r_{w},\infty)$.

\

Other vector-valued inequalities of interest are
\begin{equation}\label{eq:MR}
\Big\|
\Big(\sum_{1\le k\le N} |e^{-\zeta_{k}L} f_k|^2\Big)^{\frac12}
\Big\|_{L^q(w)}
\le
C\,\Big\|
\Big(\sum_{1\le k\le N} | f_k|^2\Big)^{\frac12}
\Big\|_{L^q(w)}
\end{equation}
for  $\zeta_{k}\in \Sigma_{\alpha}$ with $0<\alpha<\pi/2-\vartheta$  and $f_{k}\in L^p(w)$ with a constant $C$ independent of $N$, the choice of the $\zeta_{k}$'s and the $f_{k}$'s.  We restrict to $1<q<\infty$ and $w \in  A_{\infty}$ (we keep working on $\RR^n$). By a theorem of
L. Weis \cite[Theorem 4.2]{W}, we know that the existence of such a constant is  equivalent to the maximal $L^p$-regularity of $L$ on $L^q(w)$ with one/all $1<p<\infty$, that is the existence of a constant $C'$ such that for all $f\in {L^p((0,\infty), L^q(w))}$ there is  a solution $u$ of the parabolic problem on $\RR^n\times (0,\infty)$,
$$
u'(t) +Lu(t)=f(t), \ t>0, \quad u(0)=0,
$$
with
$$\|u'\|_{L^p((0,\infty), L^q(w))}  + \|Lu\|_{L^p((0,\infty), L^q(w))} \le C' \|f\|_{L^p((0,\infty), L^q(w))}.
$$

\begin{prop} Let $w\in A_\infty$ be such that $\W_w\big(p_-(L),p_+(L)\big)\neq
\emptyset$. Then for any $q\in \Int\J_{w}(L)$, \eqref{eq:MR} holds with $C=C_{q,w, L}$ independent of
$N$, $\zeta_{k}$,$f_{k}$.
\end{prop}

This result follows from an abstract result of Kalton-Weis \cite[Theorem 5.3]{KW}  together with   the bounded holomorphic functional calculus
of $L$ on  those $L^q(w)$  that we established in Theorem  \ref{theor:B-K:weights}. However, we wish to give a  different proof using extra\-polation and preceding ideas. Note that
$2$ may not be contained in $\Int\J_{w}(L)$ and the interpolation method of \cite{BK2} may not work here.

\begin{proof}  There are three steps.

\

\noindent\textit{First step: Extrapolation.} Letting $N$, $\zeta_{k}$'s and $f_{k}$'s vary at will, we denote $\F$ the family of  all   ordered pairs $(F,G)$ of the form
$$
F=\Big(\sum_{1\le k\le N} |e^{-\zeta_{k}L} f_k|^2\Big)^{\frac12}
\qquad \mbox{and}\qquad
G= \Big(\sum_{1\le k\le N} | f_k|^2\Big)^{\frac12}.
$$
Then we have for all $(F,G)\in \F$,
\begin{equation}\label{extrapol-p}
\|F\|_{L^2(u)}
\le
C_{u}\|G\|_{L^2(u)},
\qquad
\mbox{for all }u\in A_{{2}/{p_-}}\cap
RH_{\left({p_{+}}/{2}\right)'}.
\end{equation}
Recall that $2\in (p_{-}, p_{+})=\Int\J(L)$ and $u\in
A_{{2}/{p_-}}\cap RH_{\left({p_{+}}/{2}\right)'}$ means  $2\in
\W_{u}\big (p_{-}, p_{+}\big)$. In particular, $\{e^{-\zeta L}\, : \,
\zeta\in \Sigma_{\alpha}\}$ is bounded in ${ \mathcal{L}}(L^2(u))$.
This inequality is trivially checked  with $C_{u}$ equal to  the
upper bound of this family. Applying our extrapolation result
\cite[Theorem 4.7]{AM1}, we deduce that, for all $p_-<q<p_+$ and
$(F,G)\in\F$ we have
\begin{equation}\label{extrapol-q}
\|F\|_{L^q(u)} \le C_{q,u}\,\|G\|_{L^q(u)},
\qquad
\mbox{for all }u\in A_{{q}/{p_-}}\cap
RH_{\left({p_+}/{q}\right)'}.
\end{equation}
In other words,  for all $u\in A_{\infty}$ with $\W_{u}\big (p_{-},
p_{+}\big)\ne \emptyset$, \eqref{eq:MR} holds for $q\in \W_{u}\big
(p_{-}, p_{+}\big)$ with $C$ depending on $q$ and $w$. This applies
to our fixed weight $w$ of the statement with $q\in \W_w(p_-,p_+)$.
It remains to push the range of $q$'s  to all of $\Int\J_{w}(L)$.

\

\noindent\textit{Step 2: Pushing to the right.} Take $p_0 \in
\W_{w}\big (p_{-}, p_{+}\big)$,   $q_0 \in \Int\J_{w}(L)$ with
$p_0<q<q_0$.  Fix  $N$ and the $\zeta_{k}$'s. To prove \eqref{eq:MR}
for that $q$, it suffices to apply the $\ell^2$-valued version of
Theorem \ref{theor:main-w} with underlying measure $dw$ and no weight
to $\TT$ given  by
$$
\TT f= (e^{-\zeta_{1}\, L} f_{1}, \ldots, e^{-\zeta_{N}\, L} f_{N}),
\qquad  f=(f_{1}, \ldots, f_{N}),
$$
with $S=I$. To check \eqref{T:I-A} and \eqref{T:A} we use
$\A_{r}=I-(I-e^{-r^2\, L})^m$ with $m$ large enough (here, the
$\ell^2$-valued extension). Pick a ball $B$ and $f\in
(L^\infty_{c})_{\ell^2}$. Using that $e^{-t\, L} \in
\offw{p_{0}}{q_{0}}$ we can obtain \eqref{off-v-v}, replacing $\HH$
by $\ell^2$ and with $dw$ in place of $dx$. This and the fact that
$\TT \A_{r}=\A_{r}\TT $ yield \eqref{T:A}. We are left with checking
\eqref{T:I-A}. As usual we split $f$ as $\sum_{j\ge 1}
\bigchi_{C_{j}(B)}f$ (componentwise). The term with
$\bigchi_{C_{1}(B)}f= \bigchi_{4\,B}f$ is treated using the
$L^{p_{0}}_{\ell^2}(w)$ boundedness of $\TT$ (first step) and of
$I-\A_{r}$ ($\ell^2$-valued extension of Proposition
\ref{prop:sg-w:extension}). The terms $\bigchi_{C_{j}(B)}f$, $j\ge
2$, are treated using off-diagonal estimates injecting the Khintchine
inequality in the process: Let $F(x)= \| \TT
(I-\A_{r})(\bigchi_{C_{j}(B)}f)(x)\|_{\ell^2}$. Let $r_{1}, \ldots,
r_{N}$ be the $N$  first  Rademacher functions on $[0,1]$. Then, by
Khintchine's inequalities (see \cite{Ste} for instance)
\begin{align*}
F(x)
&=
\Big(
\int_{0}^1
\Big|
\sum_{1\le k\le N}  r_{k}(t)
\varphi_{\zeta_{k}}(L)(\bigchi_{C_{j}(B)}f_k)(x)
\Big| ^2\, dt
\Big)^{\frac12}
\\
&
\sim
\Big(
\int_{0}^1
\Big|
\sum_{1\le k\le N}  r_{k}(t)
\varphi_{\zeta_{k}}(L)(\bigchi_{C_{j}(B)}f_k)(x)
\Big| ^{p_{0}}\, dt
\Big)^{\frac1{p_{0}}}
\end{align*}
where $z\longmapsto \varphi_{\zeta}(z)=e^{-\zeta \, z} (1-e^{-r^2\,
z})^m$ for $\zeta\in \Sigma_{\alpha} $  is bounded on $ \Sigma_{\mu}$
when  $\vartheta <\mu <\pi/2-\alpha$. Remark that the functions
$\eta_{\pm, \zeta}$ associated to $\varphi_{\zeta}$ by
\eqref{phi-L:eta}  are easily shown to satisfy
$$
|\eta_{\pm, \zeta}(z)|
\lesssim
\frac1{|z|+|\zeta|} \min \Big (1,
\Big( \frac {r^2}{|z|+|\zeta|}\Big)^m\Big) \lesssim  \frac {\,
r^{2m}}{|z|^{m+1}},
\qquad z\in \Gamma_{\pm},
$$
where the implicit constant is  independent of $z, \zeta, r$. Thus,
using the representation \eqref{phi-L} for $\varphi_{\zeta}(L)$,
integrating $F(x)^{p_{0}}$ against $dw$ and using Minkowski's
integral inequality we obtain
$$
\Big( \aver{B} F(x)^{p_{0}} \,dw(x) \Big)^{\frac1{p_{0}}}
\lesssim
\int_{\Gamma_{+}}
\Big(
\int_0^1\aver{B} |e^{-z\, L} (\bigchi_{C_{j}(B)}h(\cdot, t,
z))(x)|^{p_{0}}\, dw(x)\,dt
\Big)^{\frac1{p_{0}}}\, |dz|
$$
 where
$$
 h(x, t, z) =  \sum_{1\le k\le N}  r_{k}(t) \,\eta_{+,\zeta_{k}}(z)\, f_{k} (x),
 $$
plus the similar term on $\Gamma_{-}$. Using $e^{-z\, L} \in
\offw{p_{0}}{p_{0}}$ for $z\in \Gamma_{+}$, the right hand side in
the above inequality is bounded by
$$
2^{j\,\theta_1} \int_{\Gamma_{+}}
\dec{\frac{2^j\,r}{\sqrt{|z|}}}^{\theta_2}\,
\expt{-\frac{\alpha\,4^j\,r^2}{|z|}}\, \Big( \int_0^1\aver{C_{j}(B)}
|h(x, t, z)|^{p_{0}}\, dw(x)\,dt
\Big)^{\frac1{p_{0}}}\, |dz|.
$$
Using again Khintchine's inequality, this is comparable  to
$$
2^{j\,\theta_1} \int_{\Gamma_{+}}
\dec{\frac{2^j\,r}{\sqrt{|z|}}}^{\theta_2}\,
\expt{-\frac{\alpha\,4^j\,r^2}{|z|}}\, \Big( \aver{C_{j}(B)} \Big(
\sum_{1\le k\le N}  | \eta_{+,\zeta_{k}}(z)\, f_{k}
(x)|^2\Big)^{\frac{p_{0}}2}\, dw(x)\Big)^{\frac1{p_{0}}}\, |dz|.
$$
At this point, we use the upper bound on $\eta_{\zeta_{k}}$ and
integrate in $z$  if $2\, m >\theta_{2}$ to obtain  that the latter
is controlled by
$$
2^{j(\theta_{1}-2\, m)}
\Big( \aver{C_{j}(B)} \Big( \sum_{1\le k\le
N}  |  f_{k} (x)|^2\Big)^{p_{0}/2}\, dw(x)
\Big)^{1/p_{0}}.
$$
The condition \eqref{T:I-A} follows readily if $2\, m >\theta_{1}$ as well.

\

\noindent\textit{Step 3: Pushing to the left}. This time, it suffices
to use the $\ell^2$-valued version of Theorem \ref{theor:SHT:small}
with underlying measure $dw$ and exponents  $p_0, q_0$ such that
$\widetilde{p}_-<q_0<\widetilde{p}_+$ and $\widehat{p}_-<p_0<q<q_0$.
Then  \eqref{SHT:small:A}  follows from the $\ell^2$-valued extension
of $e^{-t\, L} \in \offw{p_{0}}{q_{0}} $, and \eqref{SHT:small:T:I-A}
is obtained with a similar argument for the one just above to prove
\eqref{T:I-A}, by switching the role of $B$ and $C_{j}(B)$ ($j\ge
2$), and using $e^{-z\, L} \in \offw{p_{0}}{p_{0}} $.
\end{proof}

\begin{remark}
\rm When $w=1$, \eqref{eq:MR} holds for $p\in \Int\J(L)$ and recall
that this interval contains 2. Our proof contains two ways of seeing
this. First, apply the extrapolation step and specialize to $u=1$.
Second, apply steps 2 and 3 with $w=1$ and transition exponent $2$
pushing to  its right or to its left. Note that one could even reduce
things to one of those two steps by using duality as, if we denote
$\TT$ by $\TT_{L}$ then $\TT^*=\TT_{L^*}$. In \cite{BK2}, Step 3 and
duality is used. However, duality does not seem to work when $w\ne 1$
on \textit{all} of $\Int\J_{w}(L)$.
\end{remark}

\section{Commutators with bounded mean oscillation functions}\label{sec:commutators}

Let $\mu$ be a doubling measure in $\re^n$. Let $b\in \BMO(\mu)$ (BMO
is for bounded mean oscillation), that is,
$$
\|b\|_{\BMO(\mu)}
=
\sup_B \aver{B} |b-b_B| d\mu=  \sup_B \frac 1 { \mu(B)}\int_{B}
|b(y)-b_B|\, d\mu
<\infty
$$
where the supremum is taken over balls and $b_B$ stands for the
$\mu$-average of $b$ on $B$. When $d\mu=dx$ we simply write $\BMO$.
If $w\in A_\infty$ (so $dw$ is a
doubling measure) then the reverse H\"older property yields that
$\BMO(w)=\BMO$ with equivalent norms.

For $T$ a  bounded sublinear operator  in some  $L^{p_0}(\mu)$, $1\le
p_{0}\le \infty$, $b\in \BMO$,  $k\in \mathbb{N}$, we define the
$k$-th order commutator
$$
T_b^k f(x)=T\big((b(x)-b)^k\,f\big)(x),
\qquad f\in L^\infty_c(\mu),\ x\in \re^n.
$$
Note that $T_b^0=T$. If $T$ is  linear
 they can be alternatively defined by
recurrence: the first order commutator is
$$
T_b^1f(x)=[b,T]f(x)= b(x)\,T f(x)- T(b\,f)(x)
$$
and for $k\ge 2$,  the $k$-th order commutator is given by
$T_b^k=[b,T_b^{k-1}]$. As it is observed in \cite{AM1}, $T_b^k f(x)$
is well-defined almost everywhere  when $f\in L^\infty_c(\mu)$ and it
suffices to obtain boundedness with $b\in L^\infty$ with norm
depending only on $\| b\|_{\BMO(\mu)}$. We state the results for
commutators obtained in \cite{AM1}.

\begin{theor} \label{theor:comm:big}
Let $\mu$ be a doubling measure on $\re^n$, $1\le p_0<q_0\le \infty$
and $k\in \mathbb{N}$. Suppose that $T$ is a sublinear operator
bounded on $L^{p_0}(\mu)$, and let $\{\A_r\}_{r>0}$ be a family of
operators acting from $L^{\infty}_c(\mu)$ into $L^{p_0}(\mu)$. Assume
that \eqref{T:I-A} and \eqref{T:A} hold with $S=I$. Let $p_0<p<q_0$
and $w\in A_{\frac{p}{p_0}}\cap RH_{\left(\frac{q_0}{p}\right)'}$. If
$\sum_{j} g(j)\,j^k<\infty$ then there is a constant $C$ independent
of $f$ and $b\in \BMO(\mu)$ such that
\begin{equation}\label{comm-Lp:w}
\|T_b^k f\|_{L^p(w)}
\le C\, \|b\|_{\BMO(\mu)}^k\,\|f\|_{L^p(w)},
\end{equation}
for all $f\in L^\infty_c(\mu)$.
\end{theor}

\begin{theor}\label{theor:comm:small}
Let $k\in \NN$, $\mu$ be a doubling Borel measure on $\re^n$ with
doubling order $D$ and $1<p_0<q_0\le \infty$. Suppose that $T$ is a
sublinear operator and that $T$ and $T_b^m$ for $m=1,\dots,k$ are
bounded on $L^{q_0}(\mu)$. Let $\{\A_r\}_{r>0}$ be a family of
operators acting from $L^\infty_c(\mu)$ into $L^{q_0}(\mu)$. Assume
that \eqref{SHT:small:T:I-A} and \eqref{SHT:small:A} hold. If $\sum_j
g(j)\,2^{D\,j}\,j^k<\infty$, then for all $p_0<p<q_0$, there exists a
constant $C$ \textup{(}independent of $b$\textup{)} such that for all
$f \in L^\infty_c(\mu)$ and $b\in \BMO(\mu)$,
$$
\|T_b^k f\|_{L^p(\mu)} \le C\,  \|b\|_{\BMO(\mu)}^k\, \|
f\|_{L^p(\mu)}.
$$
\end{theor}

With these results in hand, we have the following theorem.

\begin{theor}\label{theor:appl:comm}
Let $w\in A_\infty$, $k\in \NN$ and $b\in \BMO$. Assume one of the
following conditions:
\begin{list}{$(\theenumi)$}{\usecounter{enumi}\leftmargin=.8cm
\labelwidth=0.7cm\itemsep=0.3cm\topsep=.3cm
\renewcommand{\theenumi}{\alph{enumi}}}

\item $T=\varphi(L)$ with $\varphi$   bounded holomorphic
on $\Sigma_\mu$,
$\W_w\big(p_-(L),p_+(L)\big)\neq\emptyset$ and $p\in\Int \J_w(L)$.

\item $T=\nabla\,L^{-1/2}$,
$\W_w\big(q_-(L),q_+(L)\big)\neq\emptyset$ and $p\in\Int \K_w(L)$.

\item $T=g_L$, $\W_w\big(p_-(L),p_+(L)\big)\neq\emptyset$ and
$p\in\Int \J_w(L)$.

\item $T=G_L$, $\W_w\big(q_-(L),q_+(L)\big)\neq\emptyset$ and
$p\in\Int \K_w(L)$.
\end{list}
Then for $ f\in L_c^\infty(\re^n)$,
$$
\|T_b^k f\|_{L^p(w)}
\le C\, \|b\|_{\BMO}^k\,\|f\|_{L^p(w)},
$$
where $C$ does not depend on $f$, $b$, and  is proportional to $\|\varphi\|_{\infty}$ in case $(a)$.
\end{theor}

Let us mention that, under kernel upper bounds assumptions,
unweighted estimates for commutators in case $(a)$ are
obtained in  \cite{DY}.

\begin{proof}[Proof of Theorem \ref{theor:appl:comm}. Part $(a)$]
We fix $p\in \Int\J_w(L)$ and take $p_0,q_0\in
\Int \J_w(L)$ so that $p_0<p<q_0$. We are going to apply Theorem \ref{theor:comm:big}
with $d\mu=dw$
and  no weight to $T=\varphi(L)$ where $\varphi$ satisfies \eqref{decay:varphi}.

   First, as $p_0\in \Int \J_{w}(L)$, Theorem
\ref{theor:B-K:weights} yields that $\varphi(L)$ is bounded on
$L^{p_0}(w)$. Then, choosing $\A_r=I-(I-e^{-r^2\,L})^m$ with $m\ge 1$
large enough, we proceed exactly as in the second case of the proof
of Theorem \ref{theor:B-K:weights}. That is, we repeat the
computations of the first case with $dw$ replacing $dx$ and using the
corresponding weighted off-diagonal estimates on balls. Applying
\eqref{eq:T:A} with $dw$ in place of $dx$ to $h=\varphi(L)$ we
conclude \eqref{T:A}. Besides, \eqref{varphi-L:w:f1} and
\eqref{eq:ref1} with $dw$ replacing $dx$ lead us to \eqref{T:I-A}
(with $S=I$). Therefore, Theorem \ref{theor:comm:big} shows the
boundedness of the commutators with $\BMO$ functions since
$\|b\|_{\BMO}(w)\approx \|b\|_{\BMO}$ as noticed earlier.

It remains to remove the assumption on $\varphi$. This is done easily
if one assumes that $b\in L^\infty$. Then  the general case  with
$b\in \BMO$ follows as mentioned above.
\end{proof}

\begin{remark}\rm
The argument is the same  as in the second
case in Theorem \ref{theor:B-K:weights} but for the whole range $\Int
\J_w(L)$ (in place of working with
$p\in(\widetilde{p}_-,\widehat{p}_+)$) since we already
proved that $\varphi(L)$ is bounded in $\Int\J_w(L)$
by Theorem \ref{theor:B-K:weights}. That is,
$T=\varphi(L)$  \textit{a posteriori}  satisfies \eqref{T:I-A} and \eqref{T:A} for $d\mu=dw$
and for {\em all} $p_0,q_0\in \Int\J_w(L)$ with $p_0<q_0$.
\end{remark}

\begin{proof}[Proof of Theorem \ref{theor:appl:comm}.  Part $(b)$]  We write
$T=\nabla L^{-1/2}$ and we already know that $T$ is bounded on $L^p(w)$ for $p\in \Int \K_{w}(L)$ by Theorem \ref{theor:ext-RT}.

First consider the case
$p\in(\widetilde{q}_-,\widehat{q}_+)$. We take $p_0, q_0$ so that
$\widetilde{q}_-<p_0<\widetilde{q}_+$ and $p_0<p<q_0<\widehat{q}_+$.
Let $\A_r=I-(I-e^{-r^2\,L})^m$ where $m\ge 1$ is an integer to be
chosen. As mentioned in the second case of the proof of Theorem
\ref{theor:ext-RT}, Lemma \ref{lemma:est-Riesz} holds
with $dw$ replacing $dx$. Thus, the hypotheses of Theorem \ref{theor:comm:big}
are fulfilled with $d\mu=dw$ and we can apply it with no weight.

Next we consider the case $p\in(\widehat{q}_-,\widetilde{q}_+)$. We
take $p_0, q_0$ so that $\widetilde{q}_-<q_0<\widetilde{q}_+$ and
$\widehat{q}_-<p_0<p<q_0$. Set $\A_r=I-(I-e^{-r^2\,L})^m$ where $m\ge
1$ is an integer to be chosen. Notice that we have just proved that
the operators $T_b^l$ for $l=0,\dots,k$ are bounded on $L^{q_0}(w)$
as  $q_0 \in(\widetilde{q}_-,\widehat{q}_+)$.  We have already
seen in the third case of the proof of Theorem \ref{theor:ext-RT}
that $T$ satisfies \eqref{SHT:small:T:I-A} and \eqref{SHT:small:A} with $d\mu=dw$.
Choosing $m$ large enough yields  the needed condition for $g(j)$
to apply Theorem \ref{theor:comm:small} with $d\mu=dw$.
\end{proof}

\begin{remark}\rm
In contrast with part $(a)$,  we do not know if
 Lemma \ref{lemma:est-Riesz} holds in the
whole range $\Int\K_w(L)$ with $dw$ replacing $dx$. Indeed,  its proof relies on an
$L^{p_0}(w)$-Poincar\'e inequality which is known only if $p_{0}>r_{w}$. We get around this obstacle with Theorem
\ref{theor:comm:small} .
\end{remark}

\begin{proof}[Proof of Theorem \ref{theor:appl:comm}. Part $(c)$] We
proceed exactly as in part $(a)$ using the arguments in Theorem
\ref{theor:square:weights}, Part $(a)$, in place of those in Theorem
\ref{theor:B-K:weights}. Details are left to the reader.
\end{proof}

\begin{proof}[Proof of Theorem \ref{theor:appl:comm}. Part $(d)$] We
follow the same scheme as in part $(b)$ using the arguments in Theorem
\ref{theor:square:weights}, Part $(b)$, in place of those in Theorem
\ref{theor:ext-RT}. Details are left to the reader.
\end{proof}

Similar results can be proved for the multilinear commutators
considered in \cite{PT} (see also \cite{AM1}) which are defined by
replacing $(b(x)-b)^k$ in $T_b^k$ by $\prod_{j=1}^k (b_j(x)-b)$ with
$b_j\in\BMO$  for $1\le j\le k$. Details are left to the reader.

\section{Real operators and power weights}

Let us  illustrate  our results on Riesz transforms in
a specific case and in particular  discuss  sharpness issues. Assume in this section that $L$ has real coefficients.
Then one knows that $q_{-}(L)=p_{-}(L)=1$, $p_{+}(L)=\infty$.

If $n=1$, one has also $q_{+}(L)=\infty$, so that we have obtained
for all $1<p<\infty$ and $w\in A_{p}$,
$$
\|L^{1/2}f\|_{L^p(w)}\sim \|f'\|_{L^p(w)}.
$$
For $p=1$, there are two weak-type (1,1) estimates for $A_{1}$
weights. In fact,
 all this can  be seen from
\cite{AT2}  where it is shown that $L^{1/2}=R\frac{d}{dx}$ and $
\frac{d}{dx}= M\tilde R L^{1/2} $ with $R$ and $\tilde R$ being
classical Calder\'on-Zygmund operators and $M$ being the operator of
pointwise multiplication by $1/a(x)$. Thus the usual weighted norm
theory for Calder\'on-Zygmund operators applies.

Let us assume next that $n\ge 2$.  In this case $q_{+}(L)>2$.  The next result will help us to study sharpness.
\begin{prop}\label{prop:sharpL}
For each $q>2$, there exists a real symmetric operator $L$ on $\RR^2$
for which $q_{+}(L)=q$.
\end{prop}

\begin{proof}
This is the example of Meyers-Kenig  \cite[p. 120]{AT1}. Let $q>2$
and set $\beta =-2/q \in (-1,0)$. Consider the operator $L=-\div A
\nabla $  obtained from $-\Delta$ by pulling back  the associated
quadratic form $\int \nabla u\cdot \nabla v$ by the quasi-conformal
application $\varphi(x)={|x|^{\beta}}{x}$,  $x \in \RR^2$. That is,
$A$ is obtained by writing out the change of variable in the relation
$$
\int A(x) \nabla u(x)\cdot \nabla v(x) \, dx= \int \nabla (u\circ \varphi^{-1})(y) \cdot \nabla (v\circ \varphi^{-1})(y)\, dy,
$$ with $u,v\in C_{0}^\infty(\re^n)$.  It is easy to see that $A$ is bounded and uniformly elliptic.
Hence, $u$ is a weak solution (in $W^{1,2}_{\rm loc}$) of $L$ if and
only if $u \circ \varphi^{-1}$ is a weak solution of $-\Delta$. In
other words, weak solutions of $L$ are harmonic functions composed
with $\varphi$. Thus, the local $L^p$ integrability of the gradient
of such a solution is exactly that of $\nabla \varphi$. The latter is
in $L^p$ near 0 if and only if $p<-2/\beta$ and is bounded locally
away from 0. Thus,   for any weak solution $u$ of $L$ defined on a
ball $2B$, $\nabla u \in L^p(B)$ for $p
<-2/\beta$ and this is optimal if $B$ is the unit ball.
With this in hand, we can apply a result by Shen \cite{Shen2}
which asserts that $q_{+}(L)$ is  the supremum of those $p$ for
which all weak solutions of $L$  defined on an arbitrary ball have
$\nabla u$ in $L^p$  locally inside that ball.  In our case,
$q_{+}(L)= -2/\beta=q$.
\end{proof}

\begin{remark}\rm
Let us also stress that if $\eta$ is a smooth compactly supported
function which is equal to 1 in a neighborhood of $0$, then
$v=\varphi \eta$ satisfies $|\nabla v(x)| \sim |x|^{\beta} $ near 0,
whereas $|L^{1/2}v(x)| \le c(1+|x|)^{-1}$. See \cite[p. 120]{AT1},
for this last fact.
\end{remark}

Let us come back to a general situation and consider the power
weights $w_{\alpha}(x)=|x|^\alpha$. Then, one has $p \in
\W_{w_{\alpha}}\big(\, 1, q_{+}(L)\big)$ if and only if
\begin{equation}
\label{eq:range} 1<p<q_{+}(L)
\qquad \mbox{and}\qquad
n\left( \frac p
{q_{+}(L)} -1 \right) < \alpha < n(p-1).
\end{equation}
For $(p,\alpha)$ tight with these relations Theorem
\ref{theor:ext-RT} yields
$$
\|\nabla f\|_{L^p(|x|^\alpha)} \lesssim
\|L^{1/2}f\|_{L^p(|x|^\alpha)}.
$$
In the latter inequality, we have in fact three parameters: $p \in
(1,\infty)$, $\alpha\in (-n,\infty)$ (for $w_{\alpha}\in A_{\infty}$) and $L$ in the family of real
elliptic operators. One can study sharpness in various ways.

Fix $L$  as in Proposition \ref{prop:sharpL} with $n=2$. The remark
following this result implies that the $L^p$ inequality can not hold
for any $(p,\alpha)$ with $-2<\alpha \le 2\big(\, \frac p {q_{+}(L)}
-1 \big)$ since in this case, one can produce an $f$($=v$) where the left
hand side is infinite and the right hand side finite.

If we fix $\alpha=0$ and $L$, then the condition $1<p<q_{+}(L)$ is
necessary (and sufficient) to obtain the $L^p$ estimate \cite{Aus}.

If we fix $p\in (1,\infty)$ and let $L$ and $\alpha$ vary, then one
can take $L=-\Delta$, in which case we are looking at the $L^p$ power
weight inequality for the usual Riesz transforms. In this case, it is
known that this forces $w_{\alpha}\in A_{p}$, hence $\alpha<n(p-1)$.

Let us consider the reverse inequalities. For a given weight $w$,
Theorem \ref{theor:reverseRiesz-w-2} says that the range of exponents
for the $L^p$ inequality contains $\W_{w}(1,\infty)$, which is the
set of $p>1$ for which $w\in A_{p}$. Hence, for $w_{\alpha}$ we have
$$
\|L^{1/2}f\|_{L^p(|x|^\alpha)}
\lesssim
\|\nabla f\|_{L^p(|x|^\alpha)}
\quad \mathrm{if} \quad -n < \alpha <
n(p-1).
$$
This is the usual range for Calder\'on-Zygmund operators. This can also
be seen from the fact proved in \cite{AT1} that $L^{1/2}=T\nabla$
where $T$ is a Calder\'on-Zygmund operator. Again for fixed $p\in
(1,\infty)$, this range of $\alpha$ is best possible by taking
$L=-\Delta$.

\end{document}